\input amstex

\input texdraw
\input epsf

\documentstyle{amsppt}
\pagewidth{5.5truein}\hcorrection{0.55in}
\pageheight{7.5truein}\vcorrection{0.75in}
\TagsOnRight
\NoRunningHeads
\catcode`\@=11
\def\logo@{}
\footline={\ifnum\pageno>1 \hfil\folio\hfil\else\hfil\fi}
\topmatter
\title Proof~of~two~conjectures~of~Ciucu~and~Krattenthaler
       on~the~enumeration~of~lozenge~tilings~of~hexagons 
       with~cut~off~corners
\endtitle
\endtopmatter
\document

\def\mysec#1{\bigskip\centerline{\bf #1}\message{ * }\nopagebreak\bigskip\par}

\def\myref#1{\item"{[{\bf #1}]}"}

\def\pf{{\it Proof.\ }}

\def\epf{\hfill{$\square$}\smallpagebreak}

\def\cite#1{\relaxnext@
  \def\nextiii@##1,##2\end@{[{\bf##1},\,##2]}%
  \in@,{#1}\ifin@\def\next{\nextiii@#1\end@}\else
  \def\next{[{\bf#1}]}\fi\next}
\def\proclaimheadfont@{\smc}

\def\pf{{\it Proof.\ }}

\define\M{\operatorname{M}}


\define\twoline#1#2{\line{\hfill{\smc #1}\hfill{\smc #2}\hfill}}
\define\twolinetwo#1#2{\line{{\smc #1}\hfill{\smc #2}}}
\define\twolinethree#1#2{\line{\phantom{poco}{\smc #1}\hfill{\smc #2}\phantom{poco}}}

\def\mypic#1{\epsffile{#1}}

\bigskip
\bigskip
\twolinethree{{\smc Mihai Ciucu\footnote{Supported in part by NSF grant DMS-1101670.}}}{{\smc Ilse Fischer\footnote{Supported in part by Austrian Science Foundation FWF, START grant Y463.}}}

\bigskip
\twolinethree{{\rm Indiana University}}{{\rm \rm Universit\"at Wien}}
\twolinethree{{\rm Department of Mathematics}}{{\rm Fakult\"at f\"ur Mathematik}}
\twolinethree{{\rm Bloomington, IN 47401, USA}}{{\rm Oskar-Morgenstern-Platz 1}}
\twolinethree{{\rm }}{{\rm A-1090 Wien, Austria}}



%
%

\define\And{1}
\define\FT{2}
\define\cutcor{3}
\define\anglepap{4}
\define\ffp{5}
\define\CLP{6}
\define\Kuo{7}
\define\Kup{8}
\define\Tri{9}
\define\MacM{10}
\define\Proctor{11}
\define\PBM{12}
\define\Sta{13}
\define\Ste{14}

\bigskip
\bigskip
{\eightpoint {\smc Abstract.} {\rm In their 2002 paper \cite{\cutcor}, Ciucu and Krattenthaler proved several product formulas for the number of lozenge tilings of various regions obtained from a centrally symmetric hexagon on the triangular lattice by removing maximal staircase regions from two non-adjacent corners. For the case when the staircases are removed from adjacent corners of the hexagon, they presented two conjectural formulas, whose proofs, as they remarked, seemed at the time ``a formidable task''. In this paper we prove those two conjectures. Our proofs proceed by first generalizing the conjectures, and then proving them by induction, using Kuo's graphical condensation method \cite{\Kuo}.}}

\bigskip
\bigskip


\define\ebaa{2.1}
\define\eba{2.2}
\define\ebb{2.3}
\define\ebc{2.4}
\define\ebd{2.5}
\define\ebe{2.6}

\define\eca{3.1}
\define\ecb{3.2}
\define\ecc{3.3}
\define\ecd{3.4}
\define\ece{3.5}
\define\ecf{3.6}
\define\ecg{3.7}
\define\ech{3.8}

\define\eda{4.1}
\define\edab{4.2}
\define\edb{4.3}
\define\edc{4.4}
\define\edd{4.5}
\define\ede{4.6}
\define\edf{4.7}
\define\edg{4.8}
\define\edh{4.9}
\define\edi{4.10}
\define\edj{4.11}
\define\edk{4.12}
\define\edl{4.13}
\define\edm{4.14}
\define\edn{4.15}
\define\edo{4.16}
\define\edp{4.17}
\define\edq{4.18}
\define\edr{4.19}
\define\edes{4.20}

\define\eea{5.1}
\define\eeb{5.2}
\define\eec{5.3}

\define\efa{6.1}
\define\efb{6.2}
\define\efc{6.3}
\define\efd{6.4}
\define\efe{6.5}
\define\eff{6.6}
\define\efg{6.7}
\define\efh{6.8}
\define\efi{6.9}


\define\tba{2.1}
\define\tbb{2.2}
\define\tbc{2.3}
\define\tbd{2.4}
\define\tbe{2.5}

\define\tca{3.1}
\define\tcb{3.2}

\define\tda{4.1}
\define\tdab{4.2}
\define\tdb{4.3}

\define\tfa{6.1}


\define\fba{2.1}
\define\fbb{2.2}
\define\fbc{2.3}

\define\fca{3.1}
\define\fcb{3.2}

\define\fda{4.1}
\define\fdb{4.2}
\define\fdc{4.3}
\define\fdd{4.4}
\define\fde{4.5}
\define\fdf{4.6}
\define\fdg{4.7}
\define\fdh{4.8}

\define\fea{5.1}

\vskip-0.05in
\mysec{1. Introduction}

In their paper \cite{\cutcor}, inspired by the classical theory of the enumeration of the symmetry classes of plane partitions (see \cite{\MacM} for the original result and \cite{\And}\cite{\Sta}\cite{\Kup}\cite{\Ste} for more recent developments), Ciucu and Krattenthaler considered the problem of counting the number tilings by unit rhombi (i.e., lozenge tilings) of regions obtained from a centrally symmetric hexagon on the triangular lattice by removing staircase-shaped regions from two of its corners. The cases when these two corners are not adjacent were solved in \cite{\cutcor}, and led to several product formulas. 

On the other hand, the case when the staircase shaped regions are removed from adjacent corners proved harder. Based on numerical data, the authors were able to guess a product formula (\cite{\cutcor,Conjecture\,A.1}) that they conjectured to give the answer to the corresponding enumeration problem. However, as the authors stated in \cite{\cutcor}, ``the fact that the result, even though given in terms of a completely explicit product, is unusually complex may indicate that proving the conjecture may be a formidable task.'' A weighted variation of this problem led to another conjectural formula, presented in \cite{\cutcor,Conjecture\,A.2}.

In the current paper we prove these two conjectures. Our approach is to first extend the conjectures to more general regions, and then prove the resulting formulas by induction, using Kuo's graphical condensation method \cite{\Kuo}.

\mysec{2. Statement of the conjectures and their generalizations}

\topinsert
\twoline{\mypic{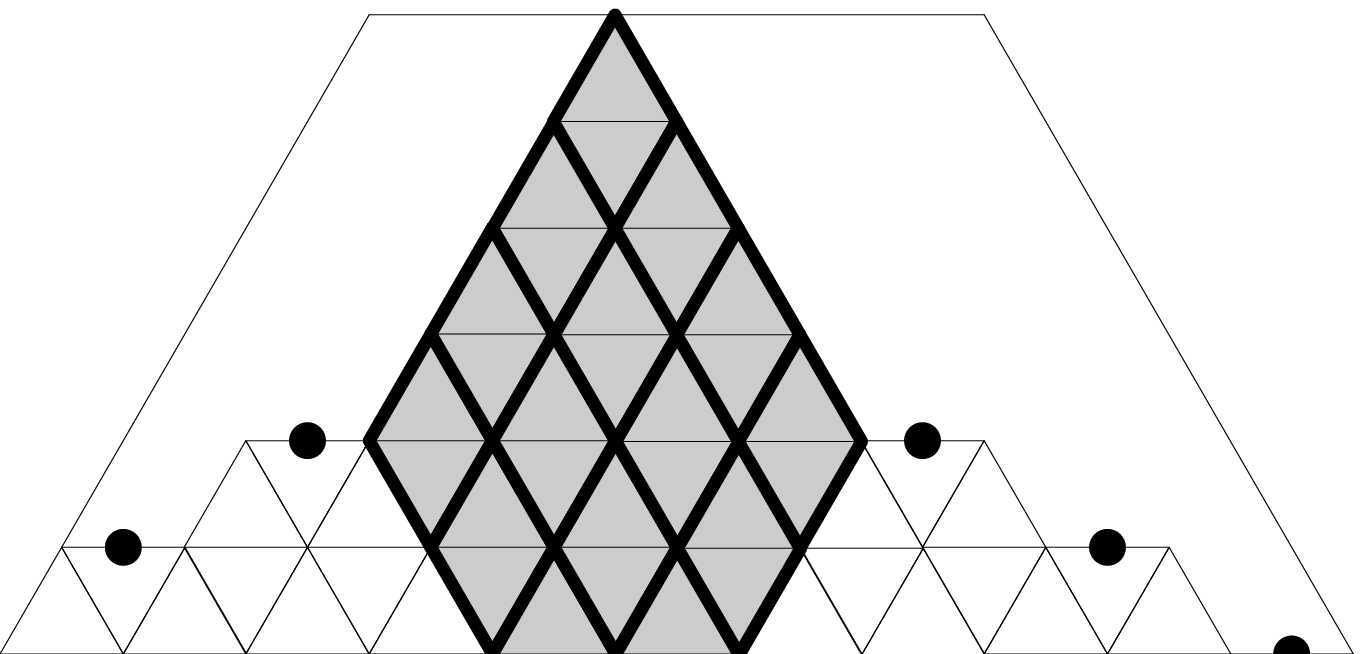}}{\mypic{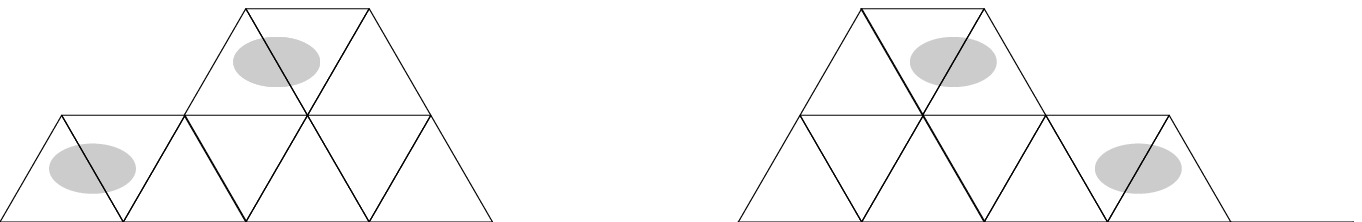}}
\twoline{Figure~\fba{\rm. \ \ \ \ \ \ }}{
Figure~\fbb{\rm. $D_{6,2,3}$}}
\vskip-5.3cm
\hbox{\hskip1.65cm
\btexdraw
  \drawdim truecm \setunitscale1.5
  \linewd.01
\htext(0.6 1){$a$}
\htext(-.75 0){$b$}
\htext(-2.1 1){$c$}
\htext(1.2 2.8){$y=(a+b-c)/2$}
\htext(4 2.8){$m=(c+b-a)/2$}
\htext(5 1){$x=a$}
\htext(3.6 -.1){$b$}
\htext(2 1){$c$}
\etexdraw
}
\vskip1cm
\endinsert

On the triangular lattice, consider the lattice hexagon of side-lengths $a$, $b$, $c$, $a$, $b$, $c$ (in cyclic order, starting from the northwestern side), and assume that $a+b+c$ is even. Let us remove the ``partial staircase'' $(a-1,a-2,\dots,(a-b+c)/2)$ from the top-left corner of this hexagon (i.e., remove northeast-to-southwest going strips consisting of $a-1$, $a-2$, etc., rhombi from that corner; see Figure {\fba}), and the partial staircase $(c-1,c-2,\dots,(c-b+a)/2)$ from the top-right corner (with the same meaning as before, except for northwest-to-southeast strips). An example is shown in Figure~{\fba}, in which the removed staircases are indicated by the white regions (the shades should be ignored at this point). We obtain a region that looks like a pentagon with an ``artificial" peak glued on top. Any lozenge tiling of this region is uniquely determined inside the rhombus that is composed out of the triangular peak and its upside-down mirror image (in Figure~{\fba} this rhombus is shaded, and the unique way to tile this rhombus is shown). Thus, for the purpose of enumerating the lozenge tilings, we may as well remove this rhombus. The leftover region now has the form of a pentagon with a notch (see Figure~{\fbb}; at this point the ellipses are without relevance). 

It will be convenient to reparametrize this region as follows. Let $x$ be the length of its southeastern edge, $y$ the number of ``bumps'' on its northwestern boundary, and $m$ the number of ``bumps'' on its northeastern boundary, and denote this region\footnote{ We point out that in \cite{\cutcor} the region obtained this way from the hexagon of side-lengths $a$, $b$, $c$, $a$, $b$, $c$ was denoted by $H_{\text n}(a,b,c)$. Thus what was denoted by $H_{\text n}(x,m+y,x+m-y)$ in \cite{\cutcor} is the same as $D_{x,y,m}$ in the current paper.}  by $D_{x,y,m}$ (note that $b$ and $c$ are then expressed in terms of $x$, $y$ and $m$ as $b=m+y$ and $c=x+m-y$). 

Note that since the depth of the notch at the top of the region $D_{x,y,m}$, which is easily seen to equal $x-y$, is non-negative, we must have $y\leq x$.

To state the conjectures we prove in this paper, we recall the definition of the Pochhammer symbol $(\alpha)_k$ (also with negative index):

$$(\alpha)_k:=\cases \alpha(\alpha+1)\cdots(\alpha+k-1)&\text {if
}k>0,\\
1&\text {if }k=0,\\
1/(\alpha-1)(\alpha-2)\cdots(\alpha+k)&\text {if }k<0.
\endcases
\tag\ebaa
$$ 
All products $\prod_{i\ge0} ^{}(f(i))_{g(i)}$ in (\eba) and (\ebb) have to be interpreted as the products over all $i\ge0$ for which $g(i)\ge0$.

For a lattice region $R$ on the triangular lattice, denote by $\M(R)$ the number of lozenge tilings of $R$.

The first conjecture we prove in this paper is the following.

\proclaim{Conjecture {\tba} (Ciucu and Krattenthaler, \cite{\cutcor,Conjecture\,A.1})}
Let $x$, $y$ and $m$ be non-negative integers with $y\leq x$.
Then the number of lozenge tilings of the region $D_{x,y,m}$ 
{\rm(}see Figure~{\fbb} for an example{\rm)}
is given by
$$\multline 
\M(D_{x,y,m})=\\
\prod _{i=1} ^{m}\frac {(x+i)!} {(x-i+m+y+1)!\,(2i-1)!}
\prod _{i=m+1} ^{m+y}\frac {(x+2m-i+1)!}
{(2m+2y-2i+1)!\,(m+x-y+i-1)!}\\
\times
{2^{\binom {m}2 + \binom y2 }}
      \prod_{i = 1}^{m-1}i!   
      \prod_{i = 1}^{y-1}i!   
      \prod_{i \ge 0}^{}
        ({ \textstyle x+i+{\frac{3}{2}} }) _{m-2i-1}    
      \prod_{i \ge 0}^{}
        ({ \textstyle   x - y+{\frac{5}{2}} + 3 i}) _{ \left \lfloor
         {\frac{3 y}{2}} -{\frac{9 i}{2}}  \right \rfloor-2}\\    
\times      \prod_{i \ge 0}^{}
        ({ \textstyle   x + {\frac{3 m}{2}} - y  + \left \lceil
         {\frac{3 i}{2}} \right \rceil+\frac {3} {2}}) 
      _{ 3 \left \lceil {\frac{y}{2}}
         \right \rceil - \left \lceil
         {\frac{9 i}{2}} \right \rceil -2}    
           \prod_{i \ge 0}
        ^{}
        ({ \textstyle {  x+ {\frac{3 m}{2}} - y + \left \lfloor
         {\frac{3 i}{2}} \right \rfloor+2}}) _{ 3 \left \lfloor
         {\frac{y}{2}} \right \rfloor - \left \lfloor {\frac{9 i}{2}} \right
         \rfloor-1}\\
\times      \prod_{i \ge 0}^{}
        ({ \textstyle  x+m - \left \lfloor {\frac{y}{2}} \right
         \rfloor}+i+1) _{  2 \left \lfloor {\frac{y}{2}} \right
         \rfloor-m - 2 i }
     \prod_{i \ge 0}
        ^{}
        ({ \textstyle  x + \left \lfloor {\frac{y}{2}} \right \rfloor+i+2})
         _{m - 2 \left \lfloor {\frac{y}{2}} \right \rfloor-2i-2} 
\\\times
{\frac{ \dsize
      \prod_{i = 0}^{y}
        ({ \textstyle x - y+3i+1}) _{m + 2 y-4i}    
      \prod_{i = 0}^{ \left \lceil {\frac{y}{2}} \right \rceil-1}
        ({ \textstyle x+m - y+i+1}) _{3 y-m-4i}    
   }
{\dsize
      \prod_{i \ge 0}^{}
        ({ \textstyle x+ {\frac{m}{2}}  - {\frac{y}{2}}+i+1}) _{y-2i}\,
  ({ \textstyle  x + {\frac{m}{2}}-
          {\frac{y}{2}}+i+{\frac{3}{2}}}) _{y-2i-1}   }}\\
\times\frac {\dsize
      \prod_{i = 0}^{y}
        ({ \textstyle x+i+2 }) _{2m - 2 i - 1}    }
 {  ({ \textstyle  x + y+2}) _{ m - y-1}  \,(m+x-y+1)_{m+y} }.
\endmultline\tag \eba
$$

\endproclaim

The second conjecture is a weighted variant of the first. Let $D'_{x,y,m}$ be the region obtained from $D_{x,y,m}$ by weighting by $1/2$ the $y$ tile positions along its nortwestern boundary and the $m$ tile positions along its norteastern boundary indicated by shaded ellipses in Figure {\fbb}. Weight each lozenge tiling of $D'_{x,y,m}$ by $1/2^k$, where $k$ is the number of lozenges in that tiling which occupy positions marked by the shaded ellipses. This weighted count of the tilings of $D'_{x,y,m}$ is the subject of the second conjecture. 

For a lattice region $R$ on the triangular lattice, some of whose lozenge positions have been weighted, we denote by $\M(R)$ the sum of the weights of its lozenge tilings (where the weight of a tiling is the product of the weights of its constituent tiles; unweighted tile positions are considered to carry weight 1).

\proclaim{Conjecture {\tbb} (Ciucu and Krattenthaler, \cite{\cutcor,Conjecture\,A.2})}
Let $x$, $y$ and $m$ be non-negative integers with $y\leq x$.
Then the weigh\-ted count of lozenge tilings of the region 
$D'_{x,y,m}$, where the lozenges along
the two zig-zag lines are weighted by $1/2$
{\rm(}see Figure~{\fbb} for an example; the lozenges that are
weighted by $1/2$ are marked by ellipses{\rm)}, 
is equal to
$$\multline 
\M(D'_{x,y,m})=\\
\prod _{i=1} ^{m}\frac {(x+i-1)!} {(x-i+m+y+1)!\,(2i-1)!}
\prod _{i=m+1} ^{m+y}\frac {(x+2m-i)!}
{(2m+2y-2i+1)!\,(m+x-y+i-1)!}\\
\times
{2^{\binom {m}2 + \binom y2 }} 
      \prod_{i = 1}^{m-1}i!   
      \prod_{i = 1}^{y-1}i!   
      \prod_{i \ge 0}^{}
        ({ \textstyle x+i+{\frac{3}{2}}}) _{m-2 i -1}    
      \prod_{i \ge 0}^{}
        ({ \textstyle x - y+3i+{\frac{7}{2}}}) _{ \left \lceil
          {\frac{3 y}{2}}-{\frac{9 i}{2}}  \right \rceil-4}\\    
\times           \prod_{i \ge 0}
        ^{}
        ({ \textstyle x+ {\frac{3 m}{2}} - y + \left \lfloor
         {\frac{3 i}{2}} \right \rfloor+\frac {3} {2}}) 
    _{ 3 \left \lceil {\frac{y}{2}}
         \right \rceil-\left \lfloor
         {\frac{9 i}{2}} \right \rfloor -1} 
      \prod_{i \ge 0}^{}
        ({ \textstyle x+ {\frac{3 m}{2}} - y + \left \lceil
         {\frac{3 i}{2}} \right \rceil+1  }) 
        _{ 3 \left \lfloor {\frac{y}{2}} \right
         \rfloor- \left \lceil
         {\frac{9 i}{2}} \right \rceil +1}  \\
\times   \prod_{i \ge 0}
        ^{}
        ({ \textstyle  x+m - \left \lfloor {\frac{y}{2}} \right
         \rfloor+i+1}) _{ 2 \left \lfloor {\frac{y}{2}} \right
         \rfloor-m - 2 i  }     \prod_{i \ge 0}
        ^{}
        ({ \textstyle  x + \left \lfloor {\frac{y}{2}} \right \rfloor+i+2})
         _{ m - 2 \left \lfloor {\frac{y}{2}} \right \rfloor-2i-2} \\
\times
{\frac{\dsize
     ({ \textstyle  x - y+{\frac{1}{2}}}) _{ \left \lfloor
      {\frac{m}{2}} \right \rfloor+2y}  ({ \textstyle x +m- y}) _{y+1}  
      \prod_{i = 0}^{y}
        ({ \textstyle x+i+1}) _{2m - 2 i }    
  }
{({ \textstyle
        x+{\frac{m}{2}} - {\frac{y}{2}}+ {\frac{1}{2}} }) 
       _{\left \lfloor {\frac{ 3y}{2}} \right \rfloor} 
       ({ \textstyle x + {\frac{3 m}{2}}  -
      {\frac{y}{2}}+1}) _{y+1}  ({ \textstyle x + {\frac{m}{2}} +
      {\frac{y}{2}}+1}) _{\left \lceil {\frac{y-2 }{2}} \right \rceil}  
}}\\
\times
\frac {  \dsize    \prod_{i = 0}^{y}
        ({ \textstyle x - y+3i+1}) _{ m + 2 y-4i}    
      \prod_{i = 0}^{ \left \lceil {\frac{y}{2}} \right \rceil-1}
        ({ \textstyle x +m- y+i+1}) _{3 y-m-4i}    } 
{  \dsize (m+x-y)_{m+y+1}\,
    ({ \textstyle x + y + \left \lceil {\frac{m}{2}} \right \rceil}) _{
     \left \lfloor {\frac{m}{2}} \right \rfloor- y + 1}  
     \prod _{i=0} ^{\lceil y/2\rceil-1}(x-y+1+3i)      }\\
\times
\frac {1}
 { \dsize     {\prod}_{i \ge 0}^{}
        ({ \textstyle x+ {\frac{m}{2}} - {\frac{y}{2}}+i+1}) _{y-2i}
  ({ \textstyle x + {\frac{m}{2}}-
          {\frac{y}{2}}+i+{\frac{3}{2}} }) _{y-2i-1}  }.
\endmultline \tag \ebb$$

\endproclaim

Our proof proceeds by first extending the above conjectures to the more general families of regions described below. Then we prove the more general formulas by induction, using Kuo's graphical condensation method described in \cite{\Kuo}.

\topinsert
\centerline{\mypic{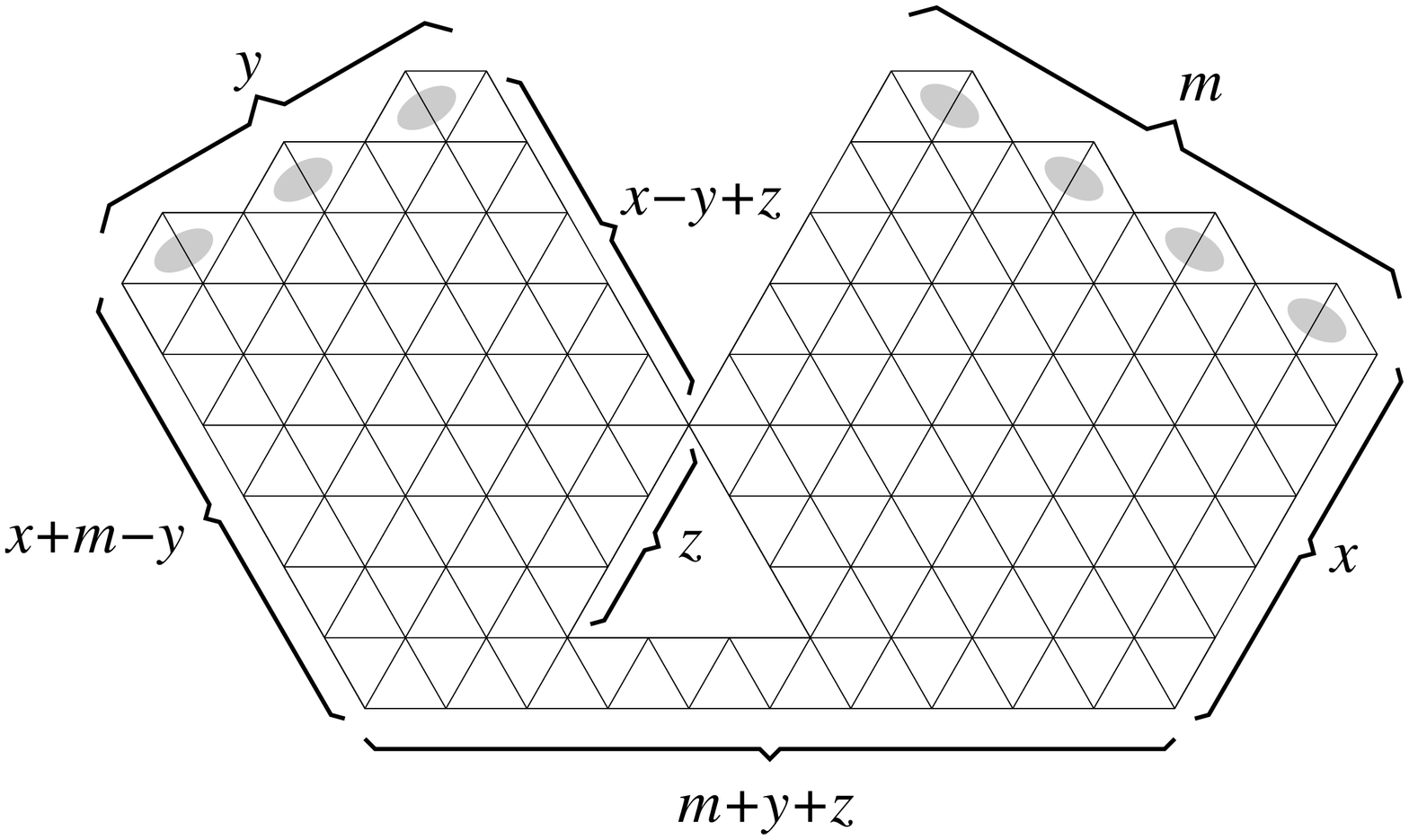}}
\centerline{{\smc Figure~{\fbc}. {\rm The generalized regions $D_{x,y,z,m}$ and $D'_{x,y,z,m}$.}}}
\endinsert

The region $D_{x,y,z,m}$ pictured in Figure {\fbc} is obtained from $D_{x,y,m}$ by introducing a triangular gap of side-length $z$ at the bottom of the notch, and modifying the side-lengths of the region as indicated in the figure (at this point the shaded ellipses should be ignored). 

Recall the encoding of lozenge tilings as families of non-intersecting paths of lozenges (see e.g. \cite{\cutcor,\S2}). Considering the paths that start from the upper left side of the triangular gap of side $z$, since they must have enough room to end at the northwestern boundary of $D_{x,y,z,m}$, we see that in order for our region to admit lozenge tilings we must have $z\leq y$. The same argument applied on the right side of the picture shows that we must also have $z\leq m$. Since the depth of the notch at the top of $D_{x,y,z,m}$ is non-negative, we must also have $x-y+z\geq0$. It then follows that $x+m-y\geq x+z-y\geq0$, so the length of the southwestern side of $D_{x,y,z,m}$ is non-negative as a consequence of the previous conditions.

The generalization of Conjecture {\tba} that we prove in this paper is the following.

\proclaim{Theorem \tbc} Let $x$, $y$, $z$ and $m$ be non-negative integers with $0\leq y-z\leq x$ and $z\leq m$. Then the number of lozenge tilings of the region $D_{x,y,z,m}$ is given by
$$
\align
&
\M(D_{x,y,z,m})=
\prod_{i=0}^{y-1}\frac{(2x+2z-i-y+1)_{y-i} i!}{(2i+1)!\,2}
\prod_{i=0}^{m-1}\frac{(2x-2y+2z+2i+2)_{m-i} (i+1)!}{(2i+2)!}
\\
&
\times
\prod_{i=0}^{m-z-1}\frac{(2x+m+3z+2i-y+3)_{y-m}(y+i-z+1)_z}{(2x-2y+3z+i+2)_{y-z}(i+1)_z}
\\
&
\times
\prod_{i=0}^{\lfloor m/3-z/3-1/3 \rfloor}(2x+y+2z+3i+3)_{3m-3z-2-9i}
\\
&
\times
\prod_{i=0}^{\lfloor m/3-z/3-2/3 \rfloor}(2x-2y+5z+6i+5)_{3m-3z-5-9i}
\\
&
\times
\prod_{i=0}^{\lfloor m/3-z/3-2/3 \rfloor}1/(2x+y-z+3m-6i-1)
\prod_{i=0}^{\lfloor m/3-z/3-1 \rfloor}1/(2x-2y+5z+6i+6)
\tag\ebc
\endalign
$$
$($if the upper index in a product is greater than the lower one, that product is taken to be $1$$)$.

\endproclaim

The extension of Conjecture {\tbb} concerns the regions $D'_{x,y,z,m}$ defined as follows. The boundary of $D'_{x,y,z,m}$ is precisely the same as the boundary of $D_{x,y,z,m}$, but the $y+m$ lozenge positions indicated in Figure {\fbc} by shaded ellipses are weighted by $1/2$. 

\proclaim{Theorem \tbd} Let $x$, $y$, $z$ and $m$ be non-negative integers with $0\leq y-z\leq x$ and $z\leq m$. Then the weighted count of the lozenge tilings of the region $D'_{x,y,z,m}$ is given by
$$
\align
&
\M(D'_{x,y,z,m})=
\prod_{i=0}^{y-1}\frac{(2x+2z-i-y)_{y-i} i!}{(2i+1)!\,2}
\prod_{i=0}^{m-1}\frac{(2x-2y+2z+2i+1)_{m-i} (i+1)!}{(2i+2)!}
\\
&
\times
\prod_{i=0}^{m-z-1}\frac{(2x+m+3z+2i-y+2)_{y-m}(y+i-z+1)_z}{(2x-2y+3z+i+1)_{y-z}(i+1)_z}
\\
&
\times
\prod_{i=0}^{\lfloor m/3-z/3-1/3 \rfloor}(2x+y+2z+3i+2)_{3m-3z-2-9i}
\\
&
\times
\prod_{i=0}^{\lfloor m/3-z/3-2/3 \rfloor}(2x-2y+5z+6i+4)_{3m-3z-5-9i}
\\
&
\times
\prod_{i=0}^{\lfloor m/3-z/3-2/3 \rfloor}1/(2x+y-z+3m-6i-2)
\prod_{i=0}^{\lfloor m/3-z/3-1 \rfloor}1/(2x-2y+5z+6i+5)
\tag\ebd
\endalign
$$
$($as in the previous theorem, if the upper index in a product is greater than the lower one, that product is taken to be $1$$)$.

\endproclaim

It is routine to check that the $z=0$ specializations of formulas (\ebc) and (\ebd) agree with formulas (\eba) and (\ebb). Thus Conjectures {\tba} and {\tbb} follow from Theorems~{\tbc}~and~{\tbd}.

\medskip
\flushpar
{\smc Remark 1.} It is interesting to note that the expression on the right hand side of formula (\ebd) is obtained from the one in (\ebc) simply by replacing $x$ by $x-\frac12$, especially since the region in the former is obtained from the region in the latter by weighting the special lozenge positions indicated in Figure {\tbc} by $\frac12$.

\medskip

Our proofs of the above results are inductive, and they rely on Kuo's powerful graphical condensation method (see \cite{\Kuo}). For ease of reference, we state below the particular instance of Kuo's general results that we need for our proofs (which is Theorem~2.4 in \cite{\Kuo}).

\proclaim{Theorem {\tbe} (Kuo)} Let $G=(V_1,V_2,E)$ be a plane bipartite graph in which $|V_1|=|V_2|+1$. Let vertices $a$, $b$, $c$ and $d$ appear cyclically on a face of $G$. If $a,b,c\in V_1$ and $d\in V_2$, then
$$
\M(G-b)\M(G-\{a,c,d\})=\M(G-a)\M(G-\{b,c,d\})+\M(G-c)\M(G-\{a,b,d\}).
\tag\ebe
$$

\endproclaim

\mysec{3. Recurrences for $\M(D_{x,y,z,m})$ and $\M(D'_{x,y,z,m})$}

As it is well known and readily seen, the lozenge tilings of a lattice region $R$ on the triangular lattice can be viewed as perfect matchings of the planar dual graph of $R$ (i.e., the graph whose vertices are the unit triangles contained in $R$, and whose edges connect unit triangles that share an edge).

Let us apply Theorem {\tbd} to the planar dual graph $G$ of the region obtained from $D_{x,y,z,m}$ by reducing the size of the triangular gap from $z$ to $z-1$, with the vertices $a,b,c,d$ chosen as indicated in Figure {\fca}. In the regions that are dual to the six graphs resulting from (\ebe), some of the lozenges are forced to be part of any tiling. What makes our proof work is the fortunate fact that after removing all the forced lozenges from these regions, each of them becomes a region of type $D$, and therefore (\ebe) turns into a recurrence for
the numbers $\M(D_{x,y,z,m})$.

\topinsert
\centerline{\mypic{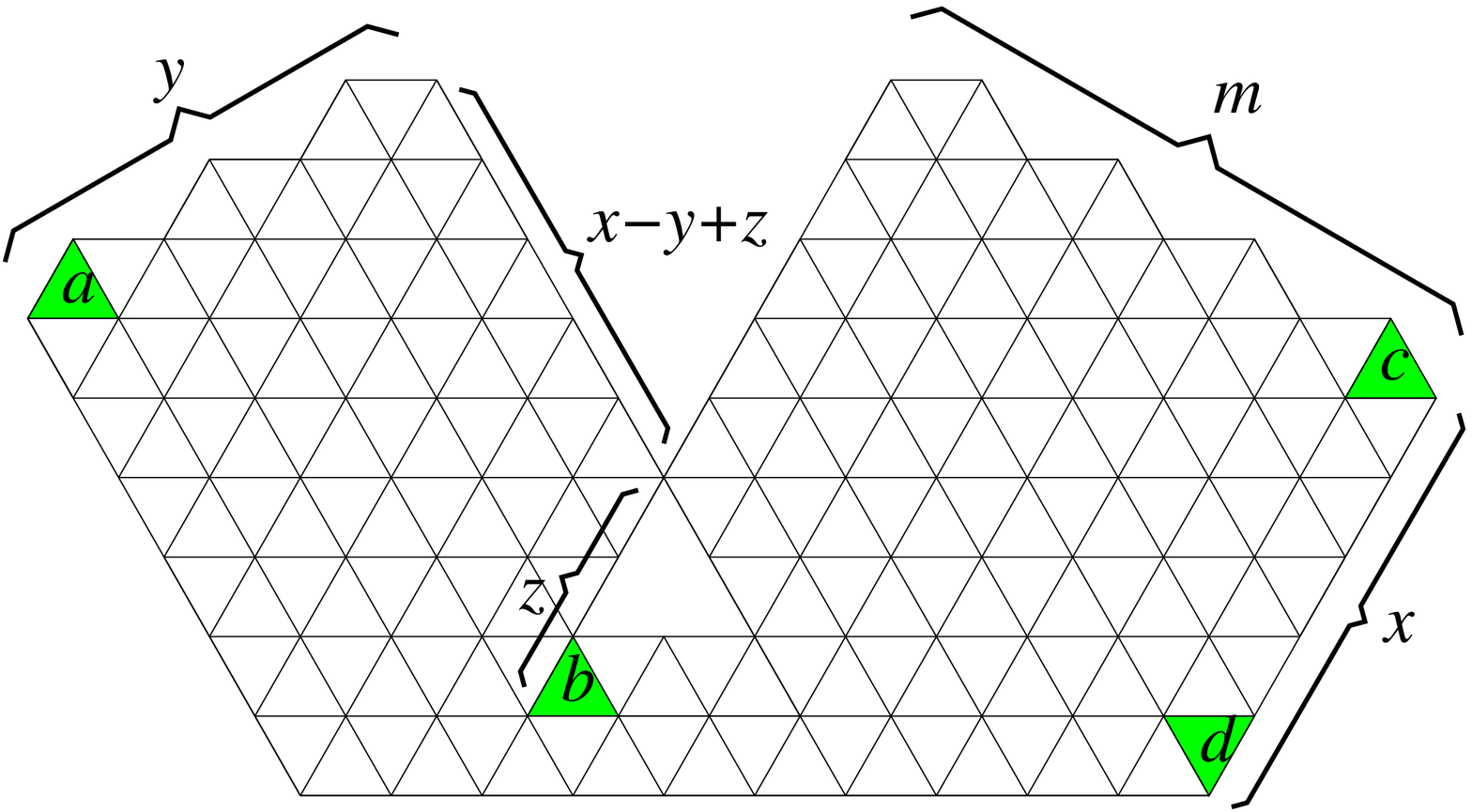}}
\medskip
\centerline{{\smc Figure~{\fca}. {\rm How we apply Kuo condensation to the regions $D_{x,y,z,m}$.}}}
\endinsert

\topinsert
\twoline{\mypic{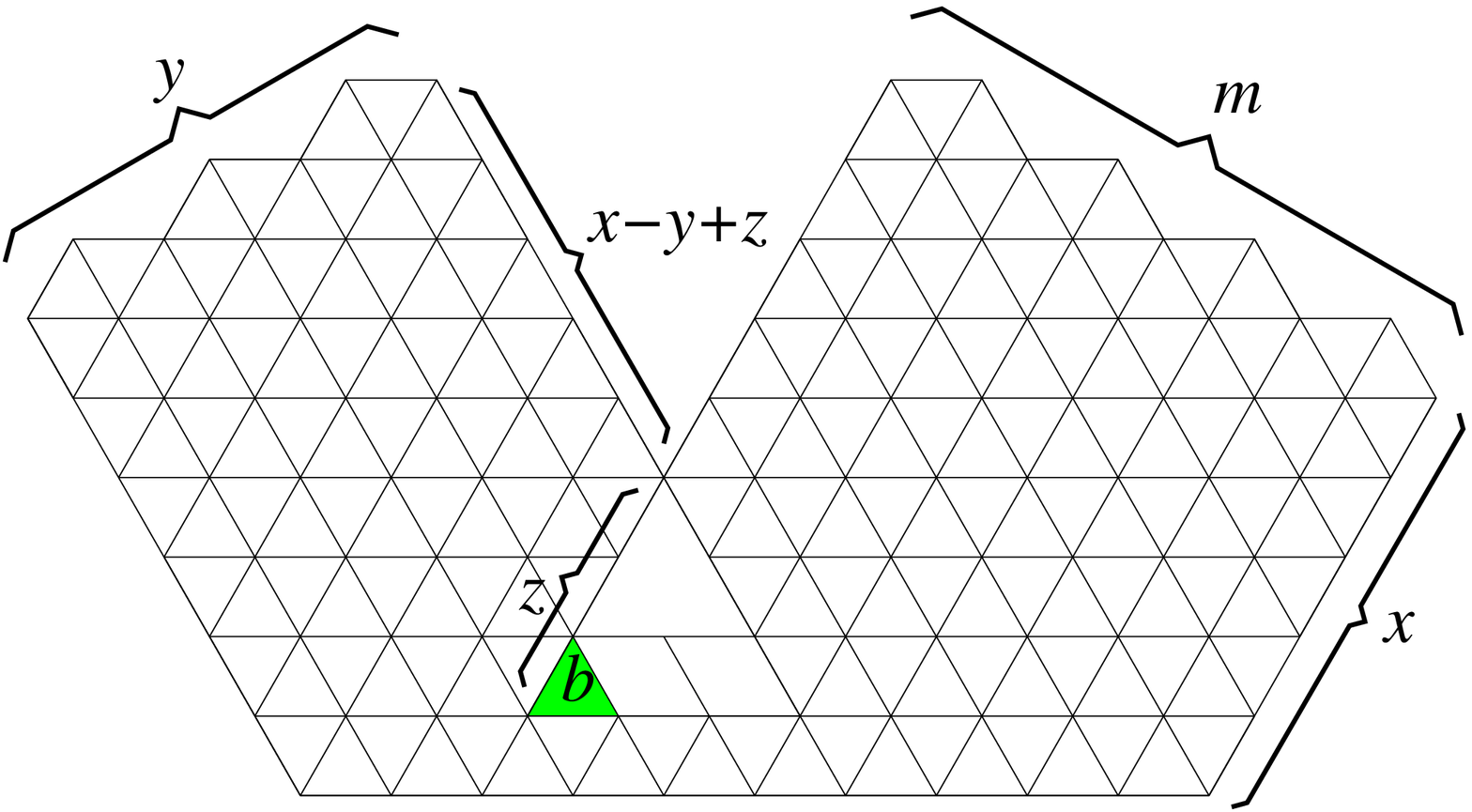}}{\mypic{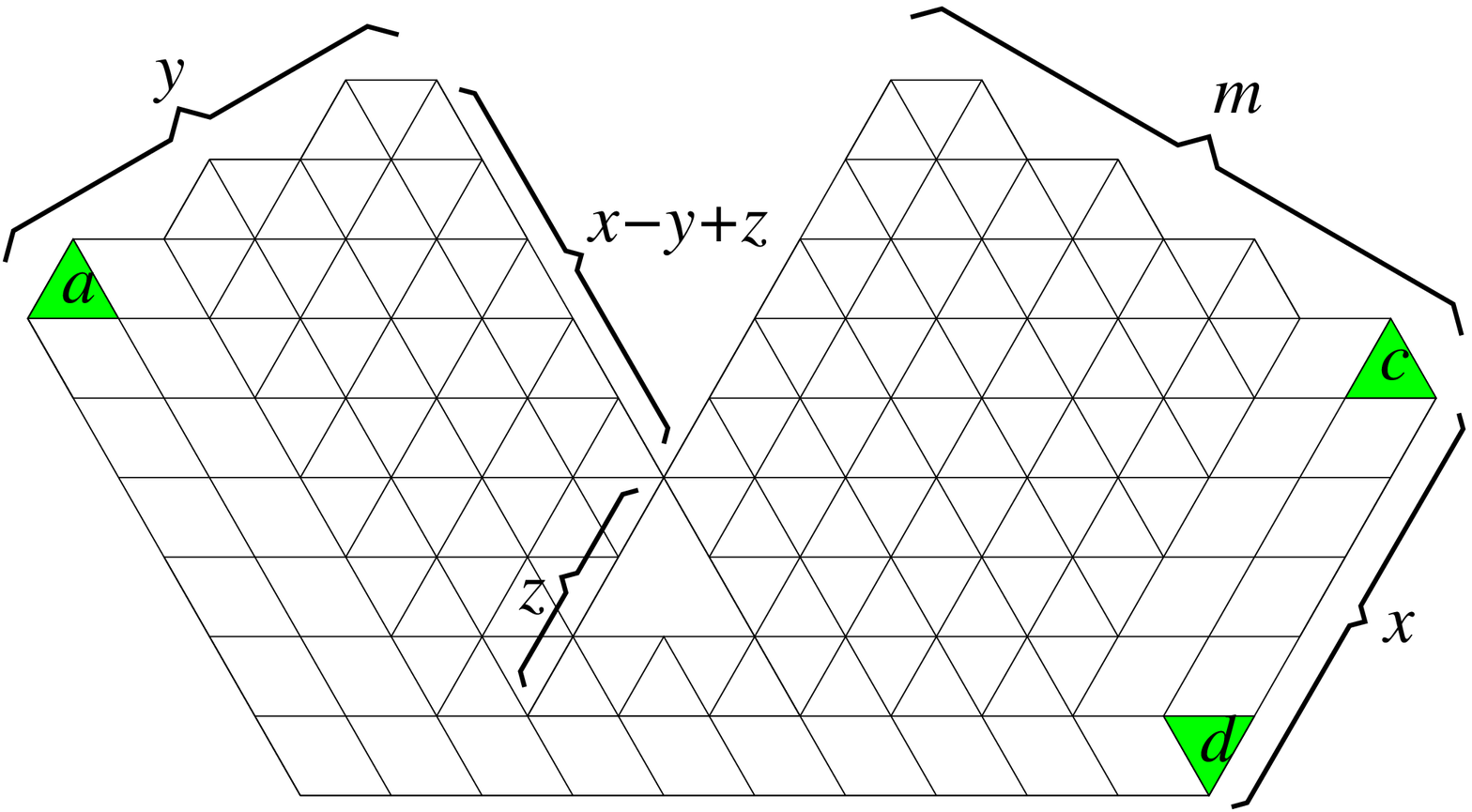}}
\twoline{\mypic{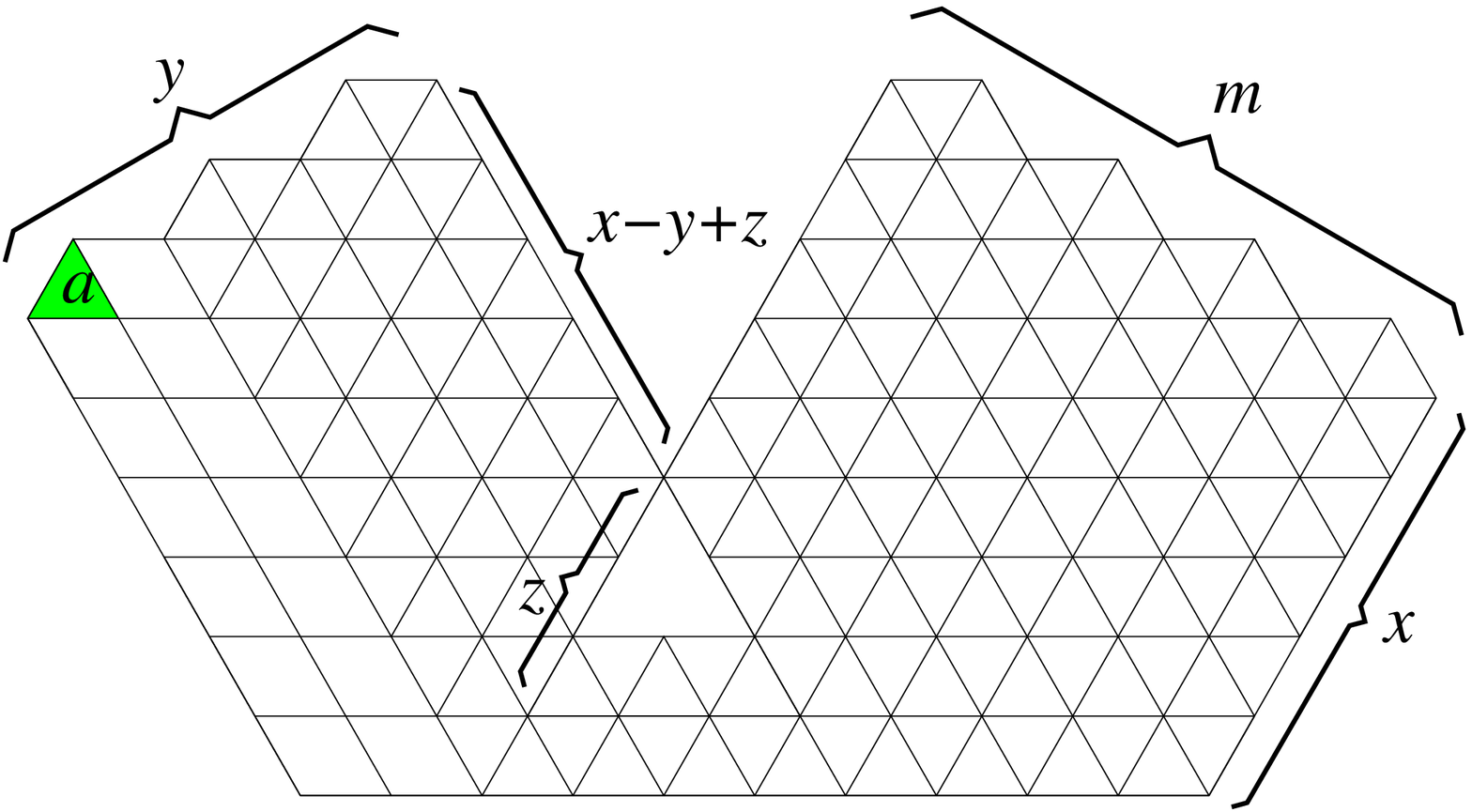}}{\mypic{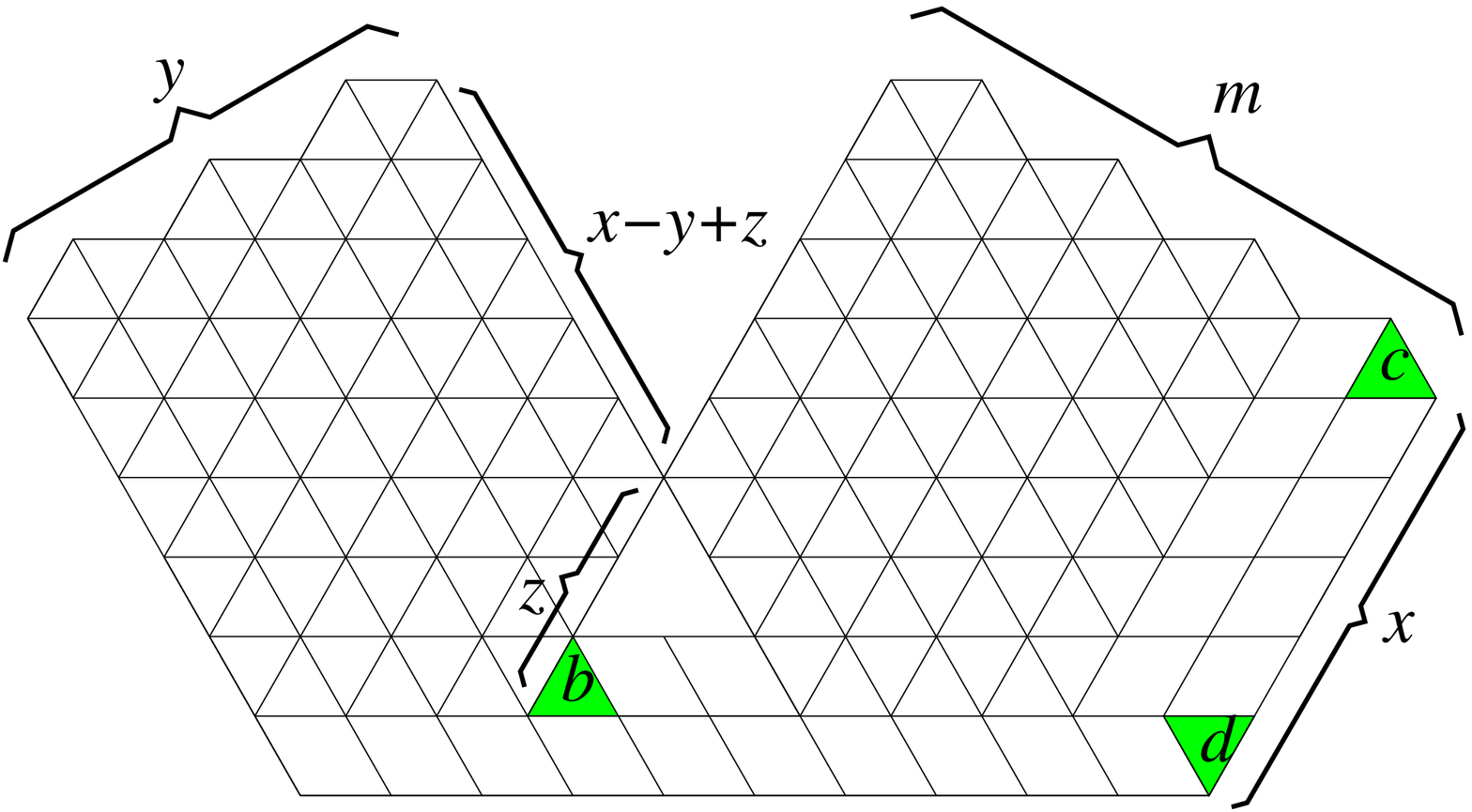}}
\twoline{\mypic{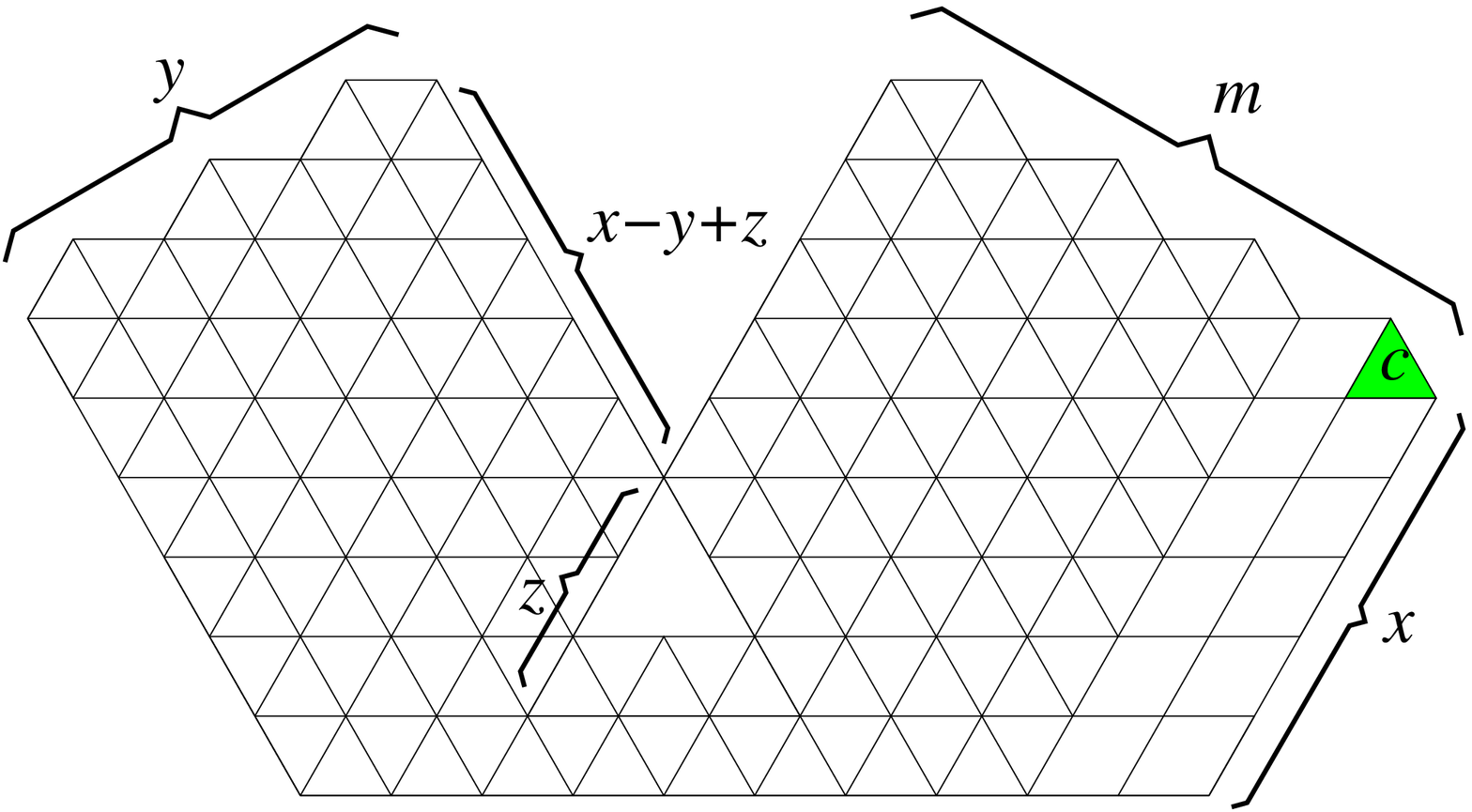}}{\mypic{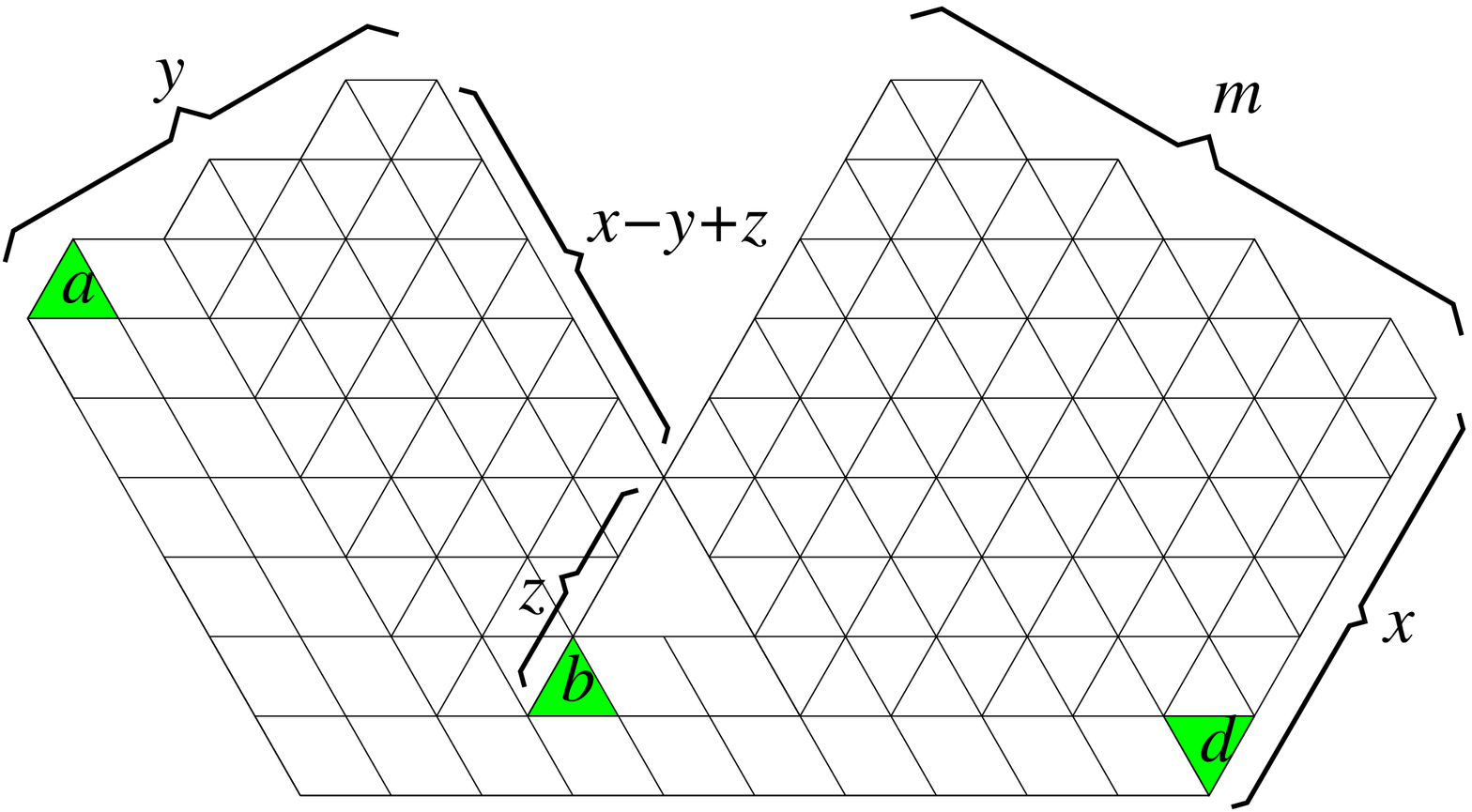}}
\medskip
\centerline{{\smc Figure~{\fcb}. {\rm Obtaining the recurrence for $\M(D_{x,y,z,m})$.}}}
\endinsert

The details are shown in Figure {\fcb}. It is apparent from it that 
$$
\M(G-b)=\M(D_{x,y,z,m}),\tag\eca
$$
$$
\M(G-\{a,c,d\})=\M(D_{x,y-1,z-1,m-1}),\tag\ecb
$$
$$
\M(G-a)=\M(D_{x,y-1,z-1,m}),\tag\ecc
$$
$$
\M(G-\{b,c,d\})=\M(D_{x,y,z,m-1}),\tag\ecd
$$
$$
\M(G-c)=\M(D_{x+1,y,z-1,m-1}),\tag\ece
$$
and
$$
\M(G-\{a,b,d\})=\M(D_{x-1,y-1,z,m}).\tag\ecf
$$

Therefore, we obtain from (\ebe) the following recurrence.

\proclaim{Lemma \tca} Let $x$, $y$, $z$ and $m$ be positive integers with $z<y$, $z<m$ and $y-z\leq x$. Then we have 
$$
\spreadlines{2\jot}
\align
&
\M(D_{x,y,z,m})\M(D_{x,y-1,z-1,m-1})
=
\\
&\ \ \ \ \ \ \ \ \ \ \ \ \ \ \ \ \ \ \ \ \ \  
\M(D_{x,y-1,z-1,m})\M(D_{x,y,z,m-1})
\\
&\ \ \ \ \ \ \ \ \ \ \ \ \ \ \ \ \ \ \   
+\M(D_{x+1,y,z-1,m-1})\M(D_{x-1,y-1,z,m}).
\tag\ecg
\endalign
$$

\endproclaim

The same approach leads to a recurrence for the numbers $\M(D'_{x,y,z,m})$. In fact, as we can see from Figure {\fcb}, none of the forced lozenges is in one of the tile positions weighted by $1/2$ in $D'_{x,y,z,m}$. It follows that the numbers $\M(D'_{x,y,z,m})$ satisfy precisely the same recurrence as the numbers $\M(D_{x,y,z,m})$.

\proclaim{Lemma \tcb} Let $x$, $y$, $z$ and $m$ be positive integers with $z<y$, $z<m$ and $y-z\leq x$. Then we have 
$$
\spreadlines{2\jot}
\align
&
\M(D'_{x,y,z,m})\M(D'_{x,y-1,z-1,m-1})
=
\\
&\ \ \ \ \ \ \ \ \ \ \ \ \ \ \ \ \ \ \ \    
\M(D'_{x,y-1,z-1,m})\M(D'_{x,y,z,m-1})
\\
&\ \ \ \ \ \ \ \ \ \ \ \ \ \ \ \ \     
+\M(D'_{x+1,y,z-1,m-1})\M(D'_{x-1,y-1,z,m}).
\tag\ech
\endalign
$$

\endproclaim

\mysec{4. The case $z=y-1$}

We start this section by recalling the following classical result due to Proctor (see \cite{\Proctor}). Let $P_{a,b,c}$ be the region obtained from the hexagon of side-lengths $a$, $b$, $c$, $a$, $b$, $c$ (clockwise starting from the northwestern side) by removing a ``maximal staircase'' from its southeastern corner (see Figure {\fda} for an illustration).

The number of lozenge tilings of $P_{a,b,c}$ is given by the following result.

\proclaim{Theorem {\tda} (Proctor \cite{\Proctor})} For any non-negative integers $a$, $b$ and $c$ with $a\leq b$ we have
$$
\M(P_{a,b,c})=
\prod_{i=1}^a
\left[
\prod_{j=1}^{b-a+1}\frac{c+i+j-1}{i+j-1}
\prod_{j=b-a+2}^{b-a+i}\frac{2c+i+j-1}{i+j-1}
\right],
\tag\eda
$$
where empty products are taken to be $1$.

\endproclaim

In fact, we will only need the special case $a=b$ of the above theorem, which we state for convenience below.

\proclaim{Corollary \tdab} For any non-negative integers $a$ and $c$ we have
$$
\M(P_{a,a,c})=
\frac{(c+1)_a}{(2c+1)_a}
\prod_{1\leq i\leq j\leq a}
\frac{2c+i+j-1}{i+j-1}.
\tag\edab
$$
\endproclaim

We present next two different generalizations of the special case $a=b$ of the above theorem, that we will need later in this section. They concern the regions $R_{x,a,k}$ and $G_{x,a,k}$ described in Figures {\fdb} and {\fdc}.

\topinsert
\centerline{\mypic{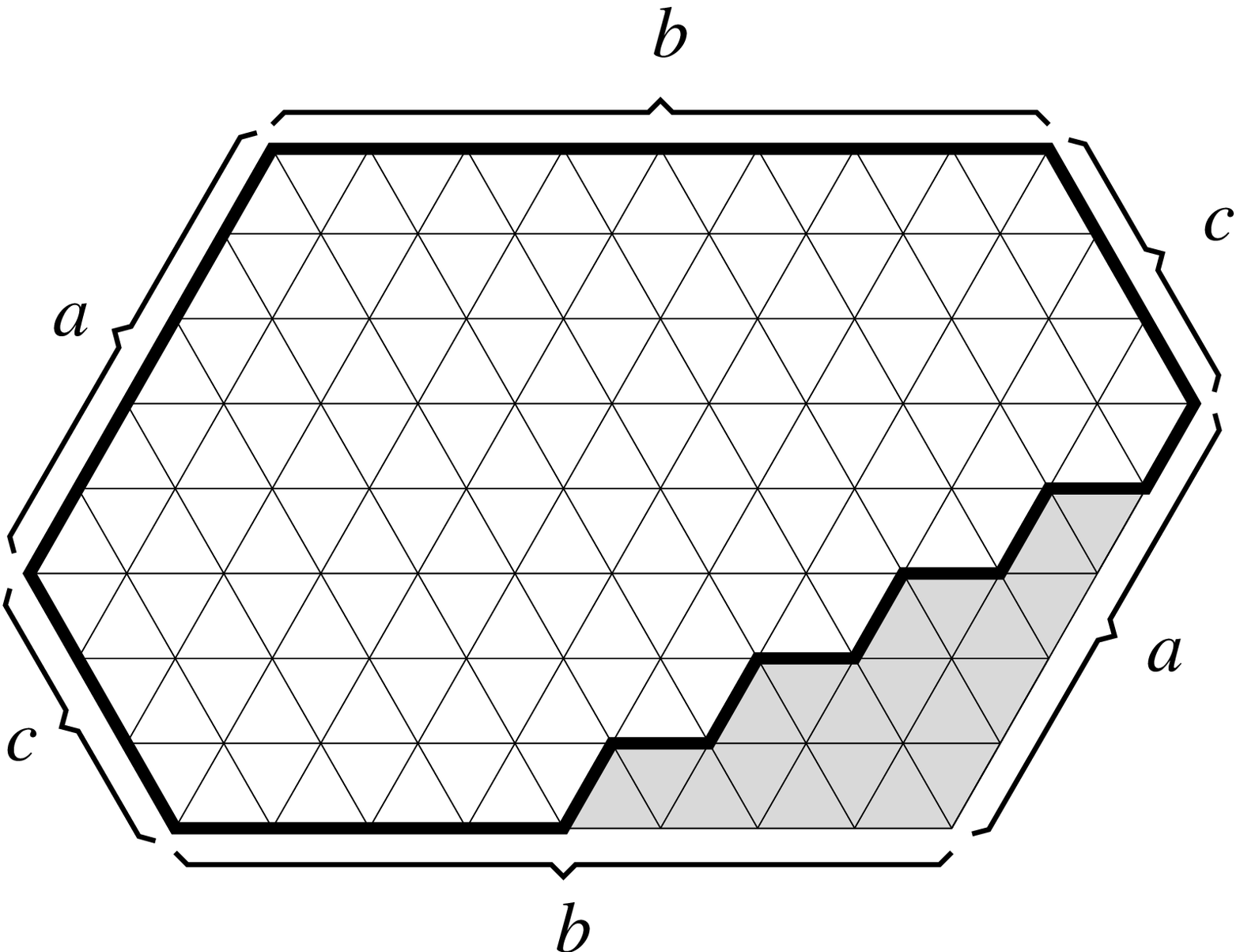}}
\medskip
\centerline{{\smc Figure~{\fda}. {\rm The region $P_{a,b,c}$ for $a=5$, $b=8$, $c=3$.}}}
\endinsert

\topinsert
\twoline{\mypic{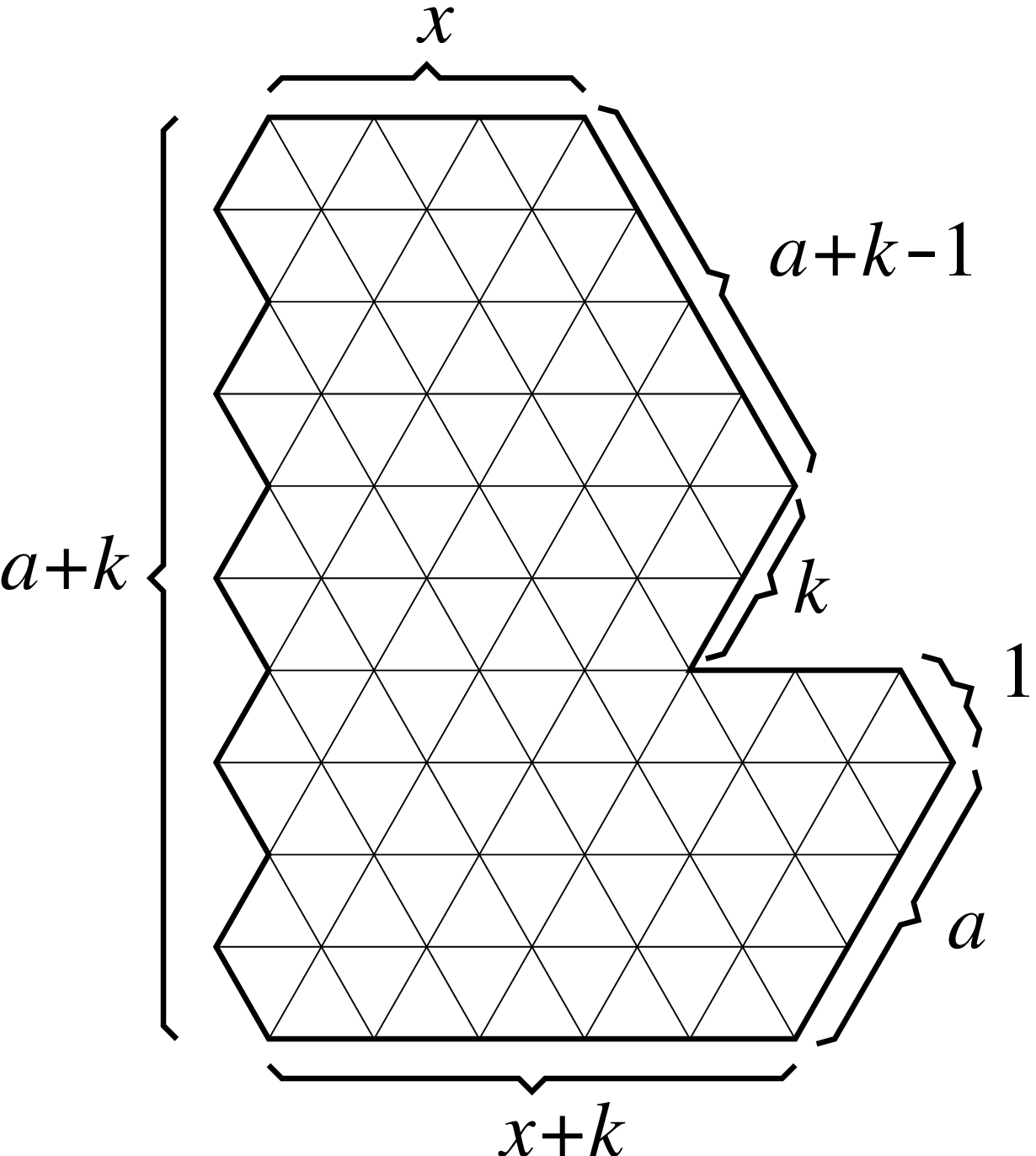}}{\mypic{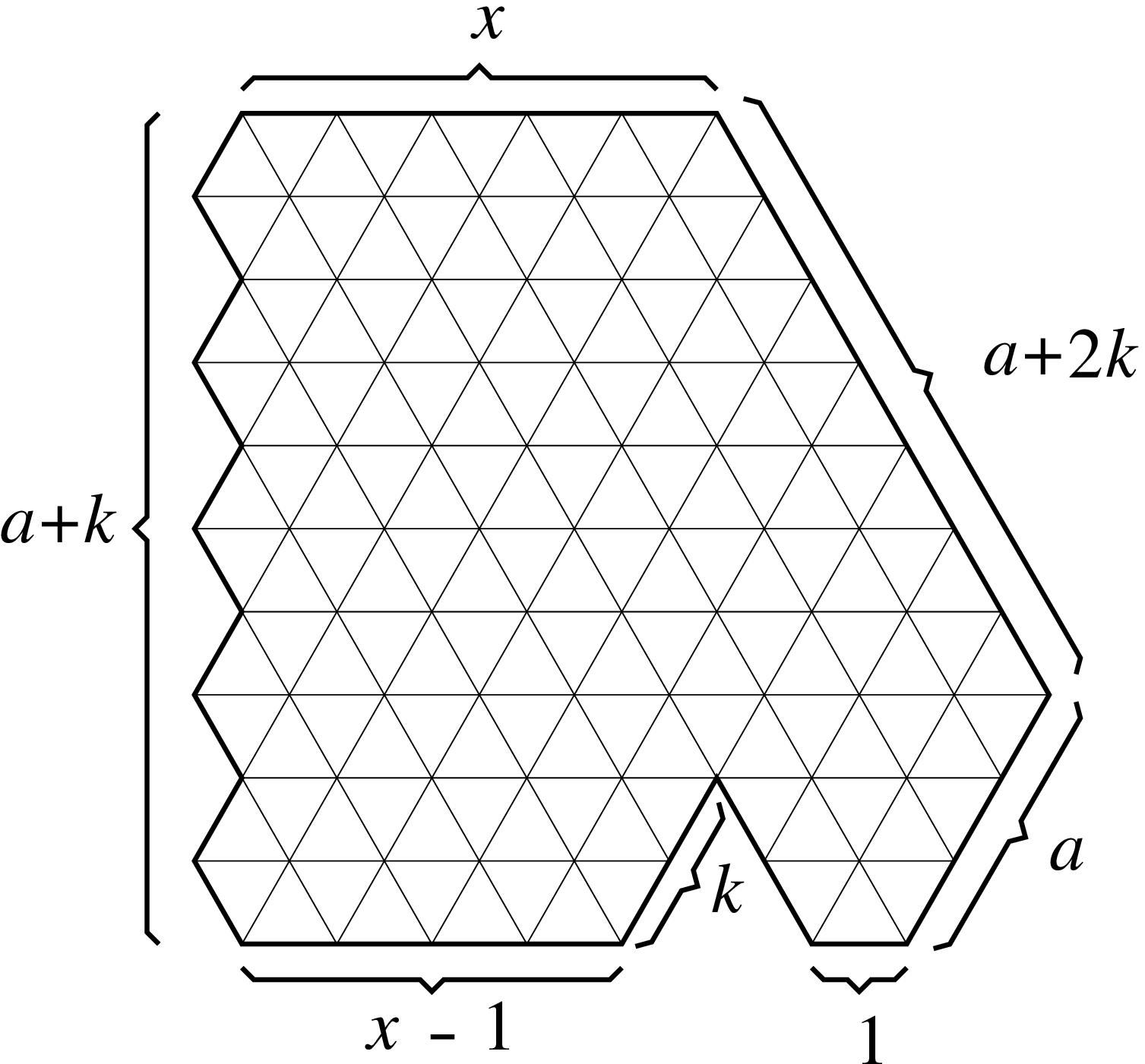}}
\medskip
\twoline{Figure~{\fdb}{. \rm $R_{x,a,k}$ for $x=3$, $a=3$, $k=2$.}}
{Figure~{\fdc}{. $G_{x,a,k}$.}}
\endinsert

\proclaim{Lemma {\tdb}} For any non-negative integers $x$, $a$ and $k$ we have

$(${\rm a}$)$.
$$
\M(R_{x,a,k})=\frac{(a+1)_k (a+k+1)_k}{(2x+a+k+1)_k k!}\M(P_{a+k,a+k,x}).
\tag\edb
$$

$(${\rm b}$)$.
$$
\M(G_{x,a,k})=\frac{(a+1)_k (2x)_k}{(2x+a+k+1)_k k!}\M(P_{a+k,a+k,x}).
\tag\edc
$$

\endproclaim

\pf (a). 
Apply graphical condensation (Theorem {\tbe}) to the region described in Figure {\fde}, with the four removed unit triangles chosen as indicated in that figure. After removing the forced lozenges in the resulting six regions, each of them turns out to be a $R_{x,a,k}$-type region; their precise parameters can be read off from Figure {\fdf}. Then (\ebe) implies that
$$
\spreadlines{3\jot}
\align
&
\M(R_{x,a,k+1})\M(R_{x+1,a+k-1,0})
=
\\
&\ \ \ \ \ \ \ \ \ \ \ \ \ \ 
\M(R_{x+1,a,k})\M(R_{x,a+k,0})
+
\M(R_{x,a+k+1,0})\M(R_{x+1,a-2,k+1}),
\tag\edd
\endalign
$$
for all non-negative integers $x$, $a$ and $k$ with $a\geq2$ (so that all regions in the above equations are defined).

Note that in fact (\edd) holds also for $a=1$, provided we take 
$$
\M(R_{x+1,-1,k+1})=0.
\tag\ede
$$
Indeed, (\ebe) can be applied as indicated in Figure {\fde} also for $a=1$, but then the region on the bottom right in Figure {\fdf} --- the entire region, before removing any forced lozenges --- has no lozenge tilings.

\topinsert
\centerline{\mypic{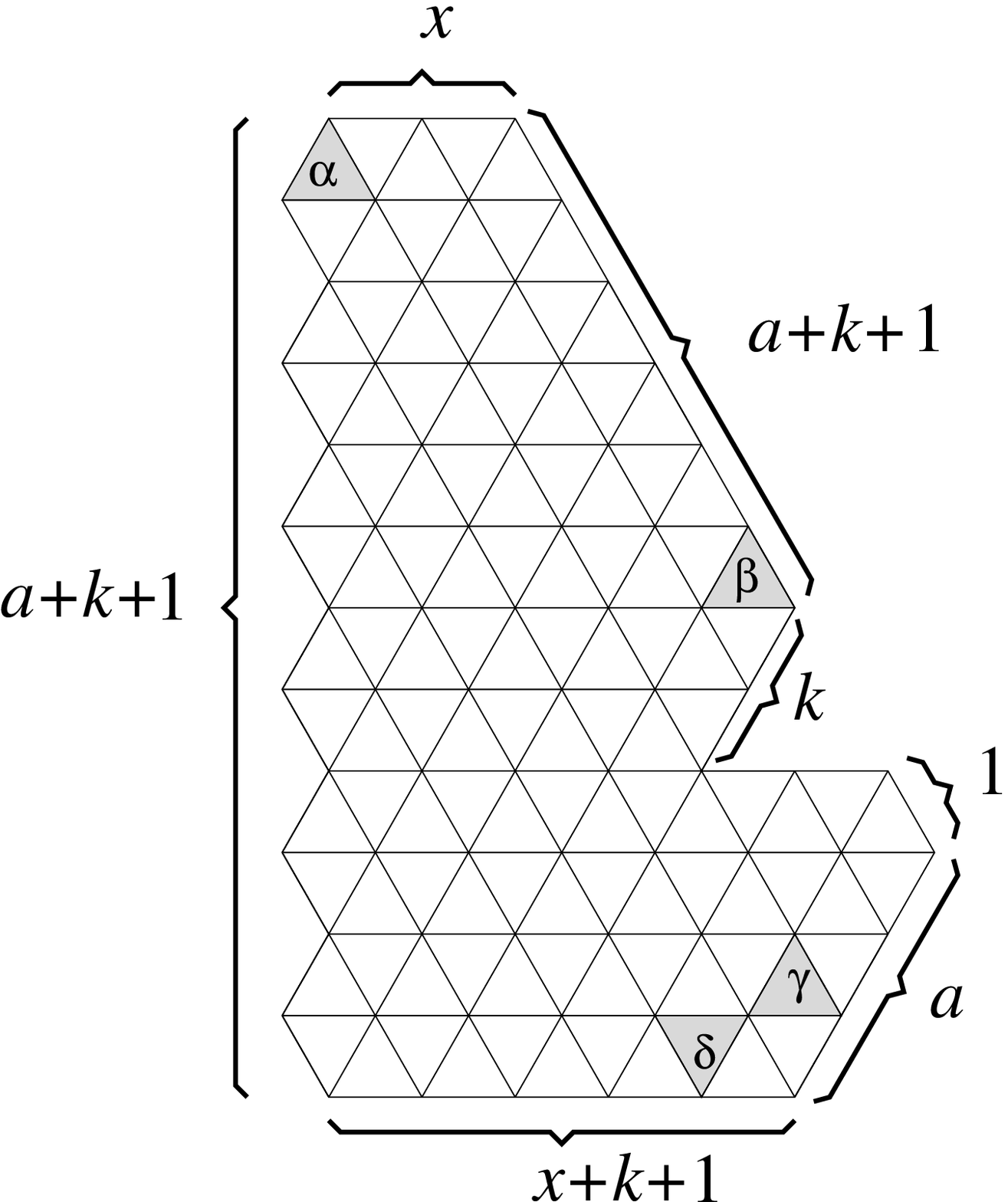}}
\medskip
\centerline{{\smc Figure~{\fde}. {\rm Region to which we apply graphical condensation.}}}
\endinsert

\topinsert
\twoline{\mypic{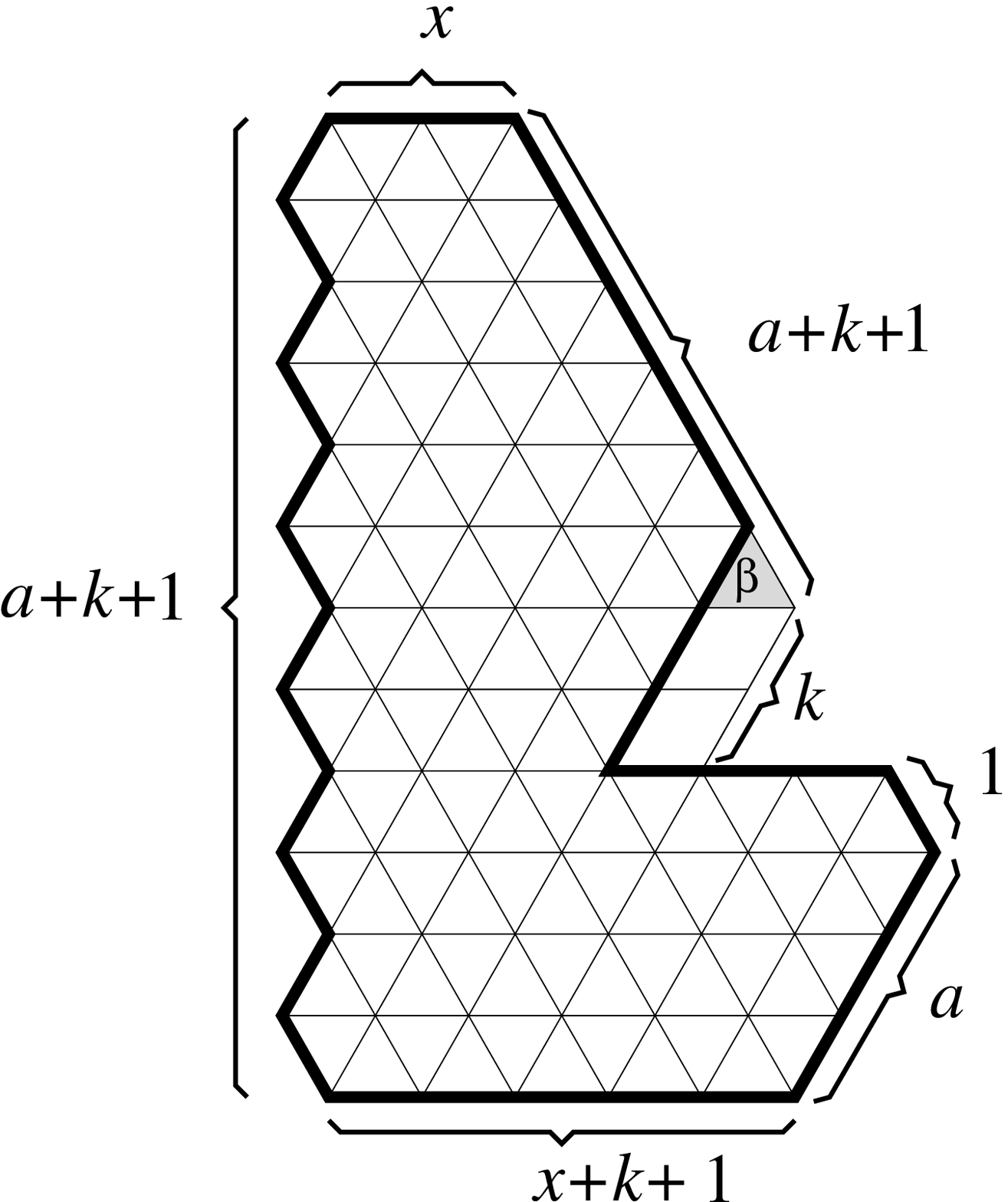}}{\mypic{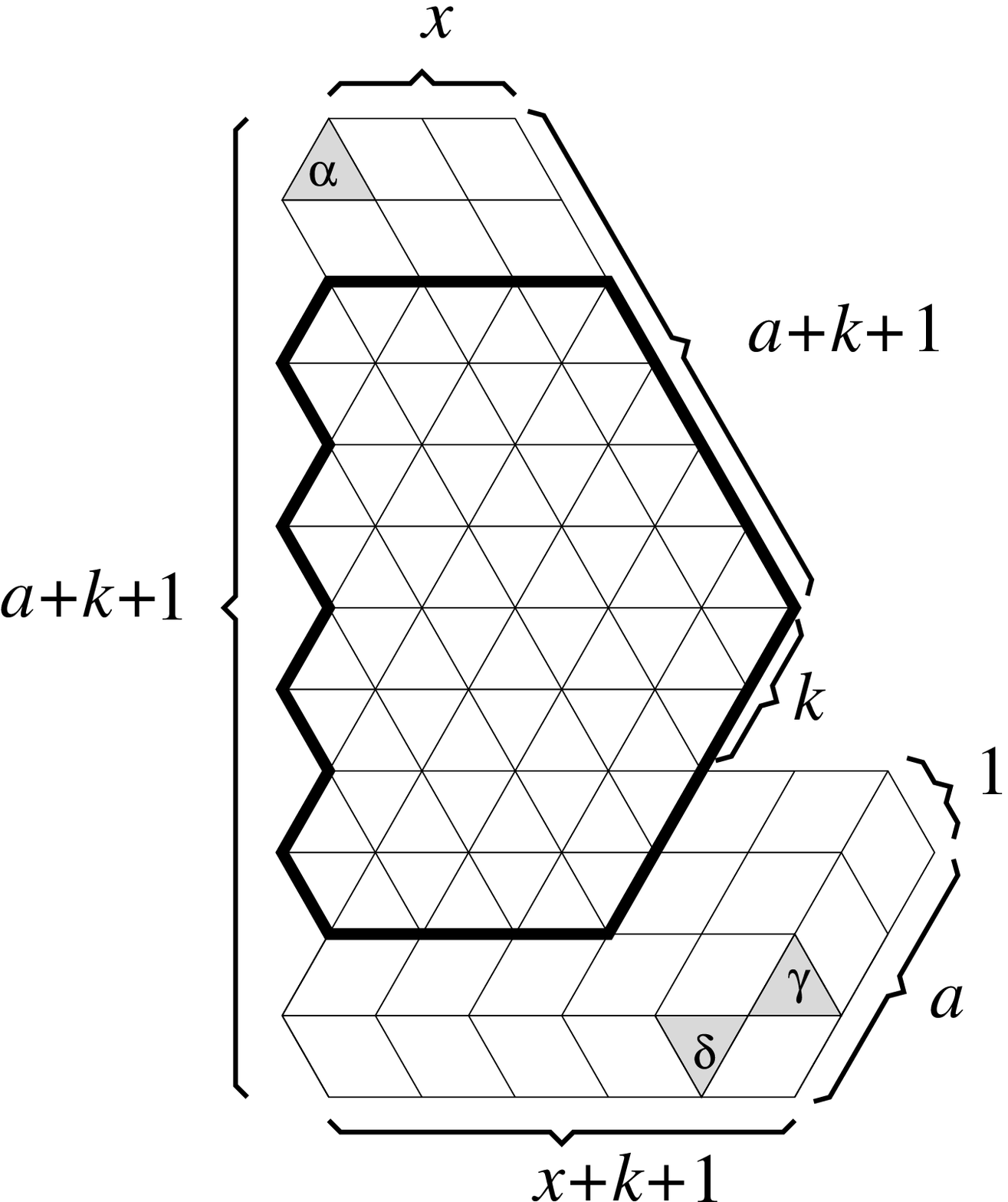}}
\twoline{\mypic{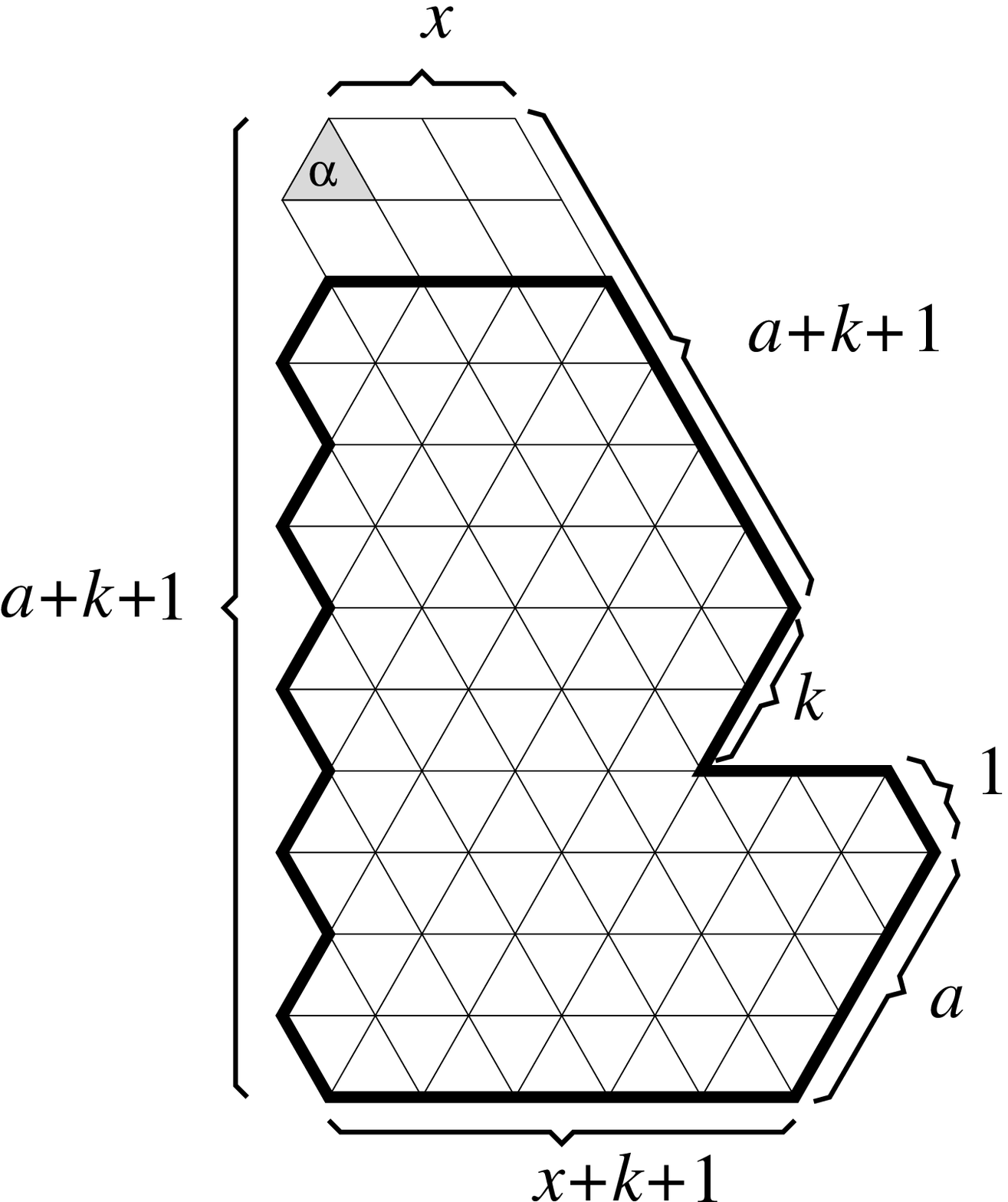}}{\mypic{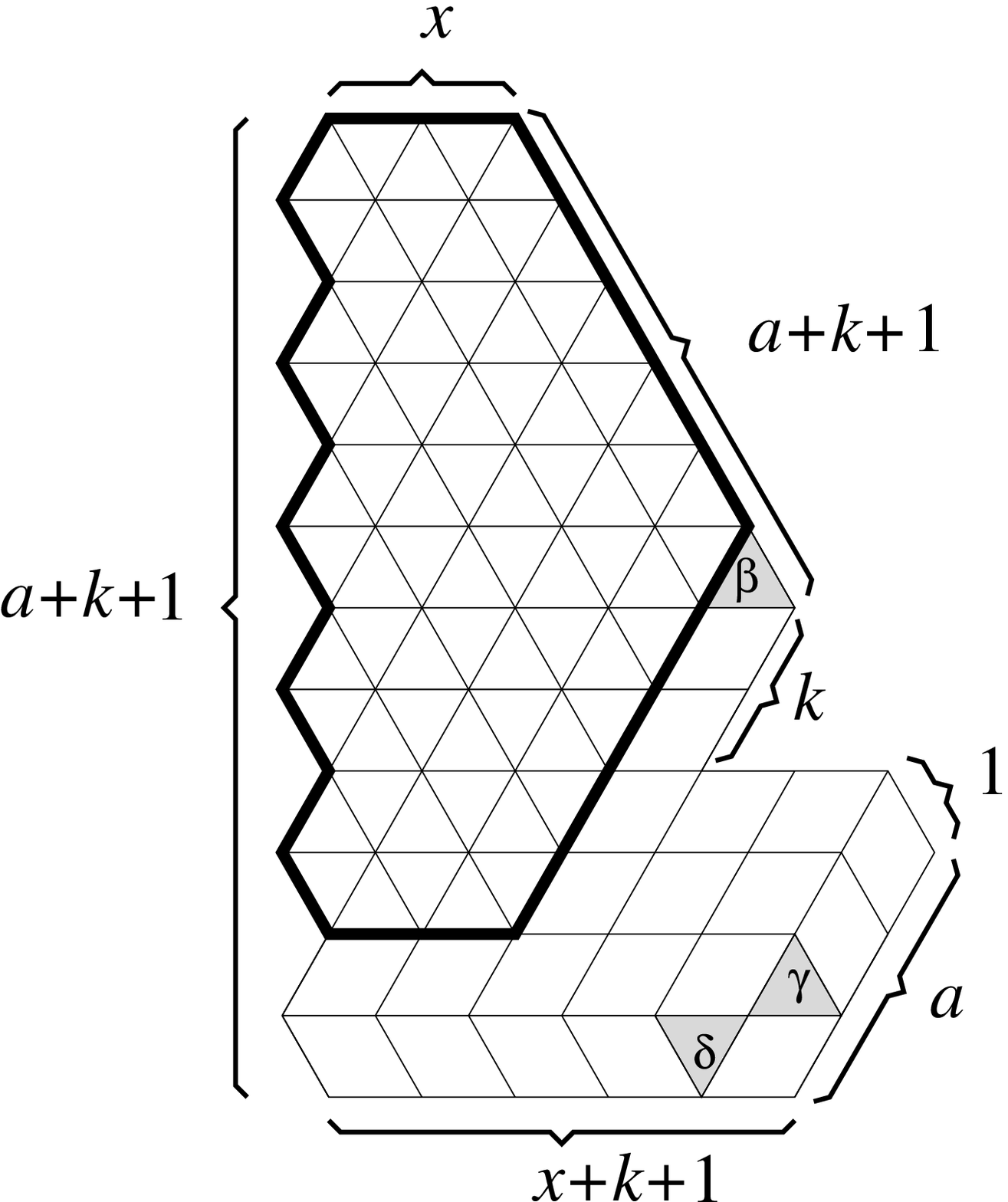}}
\twoline{\mypic{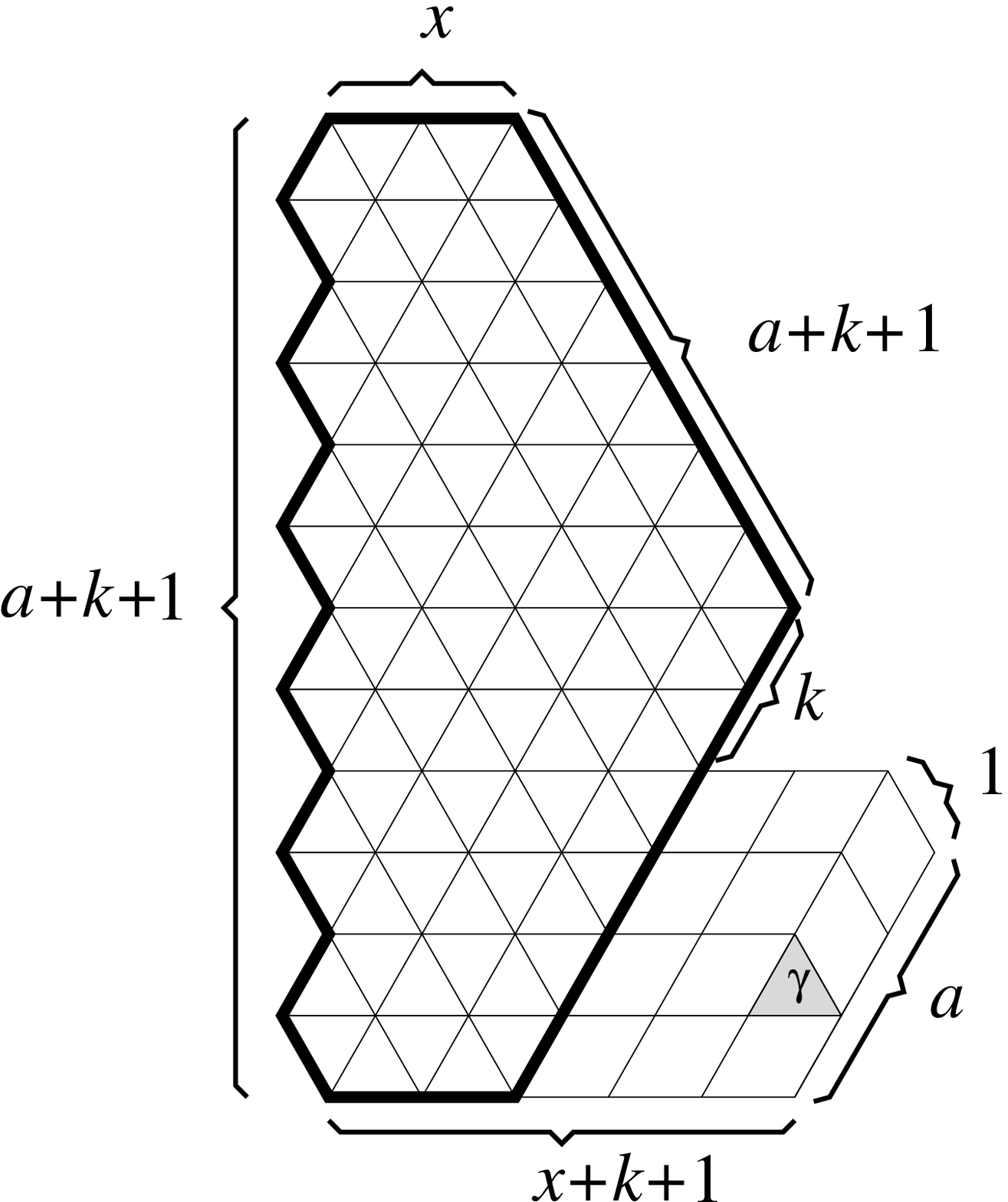}}{\mypic{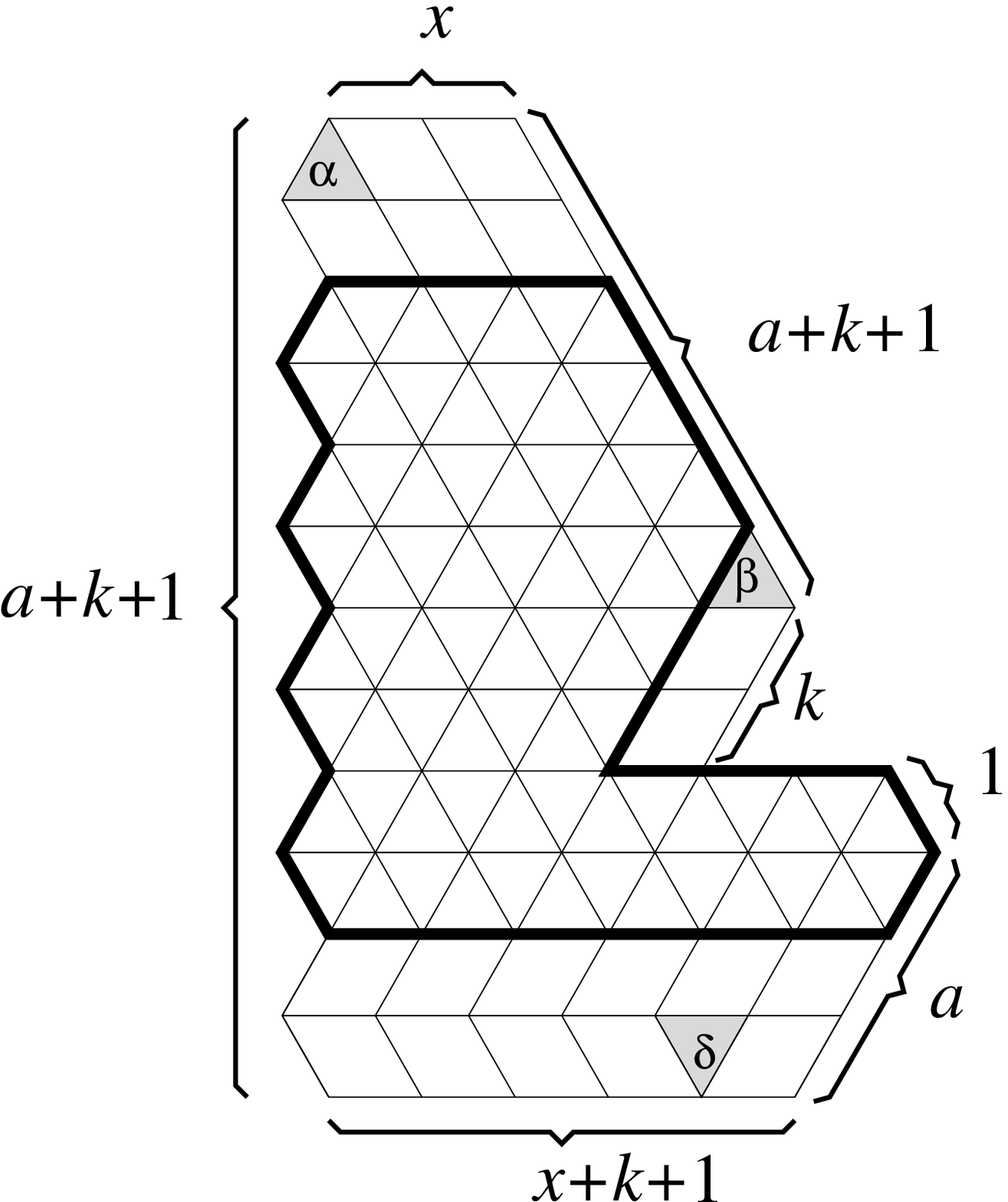}}
\medskip
\centerline{{\smc Figure~{\fdf}. {\rm Obtaining the recurrence for the $\M(R_{x,a,k})$'s.}}}
\endinsert

We use the recurrence (\edd) to prove (\edb) by induction on $a+k$. If $a+k=0$, both regions in (\edb) are empty, and all the factors in the fraction are equal to 1, so (\edb) holds. If
$a+k=1$, then either $a=1,k=0$ or $a=0,k=1$. In the former case, the fraction in (\edb) equals 1 again, and since the region $R_{x,1,0}$ is the same as $P_{1,1,x}$, (\edb) checks. In the latter, the fraction is $2/(2x+2)$, $\M(P_{1,1,x})=x+1$, and $\M(R_{x,0,1})=1$, so (\edb) holds.

For the induction step, let $n>0$ be a positive integer and assume that (\edb) holds for all non-negative $x,a,k$ with $a+k\leq n$. We need to show that (\edb) holds also for any non-negative $x,a,k$ with $a+k=n+1$.

Consider an $R$-region so that the sum of the $a$- and $k$-parameters is $n+1$. 
Suppose $a=0$. Observing that the region obtained from $R_{x,0,k}$ after removing all forced lozenges is the same as $P_{k-1,k-1,x}$, equation (\edb) is readily verified using (\eda). Similarly, if $k=0$, then the region $R_{x,a,k}$ becomes the region $P_{a,a,x}$, and (\edb) follows again by (\eda).

We can therefore assume that $a\geq1$ and $k\geq1$. Then (\edd) implies
$$
\align
\M(R_{x,a,k+1})=
\frac{\M(R_{x,a+k,0})}{\M(R_{x+1,a+k-1,0})}\M(R_{x+1,a,k})
+
\frac{\M(R_{x,a+k+1,0})}{\M(R_{x+1,a+k-1,0})}\M(R_{x+1,a-2,k+1}).
\tag\edf
\endalign
$$
Note that as the $k$-parameter of our $R$-region is at least 1, it can be written in the form of the left hand side of (\edf), with $k\geq0$. 
The fractions on the right hand side of (\edf) have explicit expressions by Corollary {\tdab}. Furthermore, since the sum of the $a$- and $k$-parameters for the remaining two regions is one or two units less than for the region on the left hand side, the induction hypothesis implies that the 
number of their lozenge tilings is given by the formula (\edb).
This way (\edf) leads to an explicit expression for $\M(R_{x,a,k+1})$, which is readily verified to agree with the expression obtained from the right hand side of (\edb) by replacing $k$ by $k+1$. This concludes the induction step, and therefore the proof of (\edb).

(b). Formula (\edc) can be proved by a perfectly similar argument. We note that it also follows as a special case of a formula of Lai (see \cite{\Tri}). \epf

\topinsert
\centerline{\mypic{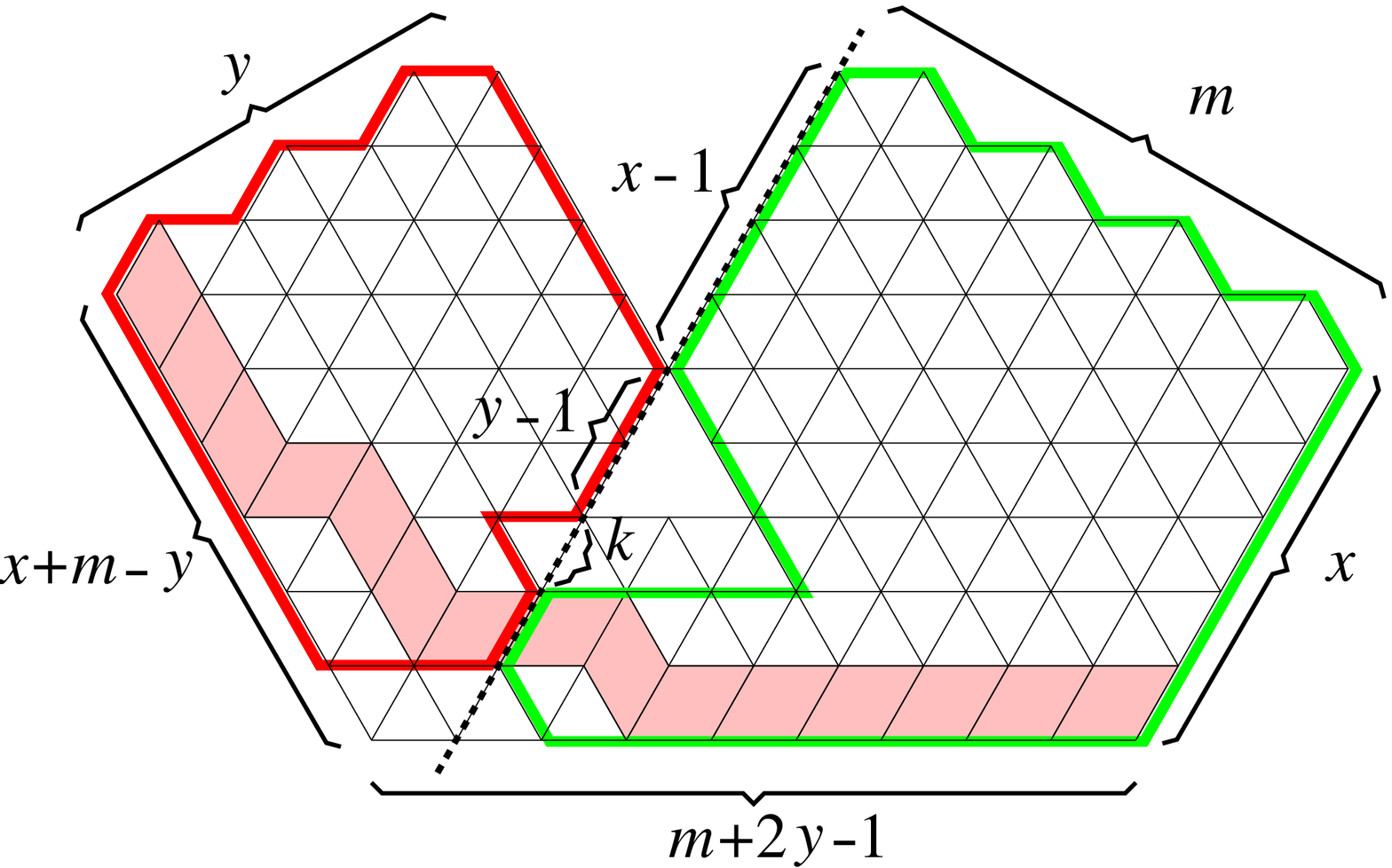}}
\medskip
\centerline{{\smc Figure~{\fdd}. {\rm Reduction to smaller regions.}}}
\endinsert

{\it Proof of the case $z=y-1$ of Theorem {\tbc}.} Set $z=y-1$, and consider the region $D_{x,y,y-1,m}$. The lengths of its sides are shown in Figure {\fdd}. 


We present the details in the case when the extension $\ell$ of the northwestern side of the triangular gap (see the dotted line in Figure {\fdd}) leaves the region through its base (and not through an interior point of its southwestern side). The only modification in the other case is that the range for $k$ in (\edg) and (\edh) is in that case between $0$ and $y$. However, due to the presence of the factor $(y-k+1)_k$ at the numerator of the summand in (\edh), the upper limit in the sum can then be replaced by infinity without changing its value, and the argument follows by exactly the same calculations as in the case we are presenting. 

Recall that each tiling of $D_{x,y,y-1,m}$ can be encoded as a family of $x$ non-intersecting paths of lozenges, starting at the unit segments on the southeastern edge. Since the southeastern side of the notch has length $x-1$, precisely one of these paths ends at the northwestern zig-zag boundary of $D_{x,y,y-1,m}$; in Figure {\fdd}, this path is indicated by shaded lozenges. 

Note that the number of lozenge tilings of $D_{x,y,y-1,m}$ for which this path of lozenges crosses $\ell$ at a fixed unit segment is equal to the product of the number of tilings of the two regions shown in thick contours just above and just below $\ell$. Using the indicated lengths in Figure {\fdd}, one sees that the region above $\ell$ is precisely the earlier defined region $R_{x-1,y-k,k}$, while the region below $\ell$ is $G_{x,m-(k+(y-1)),k+(y-1)}$.

By Figure {\fdd}, the length of the portion of $\ell$ below the triangular gap and inside the region is $x+m-(x-1)-(y-1)=m-y+2$, so the range for the parameter $k$ in Figure {\fdd} is between 0 and $m-y+1$. 

We obtain therefore that
$$
\M(D_{x,y,y-1,m})=\sum_{k=0}^{m-y+1}\M(R_{x-1,y-k,k})\M(G_{x,m-k-y+1,k+y-1}).
\tag\edg
$$
Using Lemma {\tdb}, this implies
$$
\spreadlines{3\jot}
\align
&
\M(D_{x,y,y-1,m})=
\M(P_{y,y,x-1})\M(P_{m,m,x})
\\
&\ \ \ \ \ \ \ \ \ \ \ \ 
\times
\sum_{k=0}^{m-y+1}
\frac{(y-k+1)_k(y+1)_k}{(2x+y-1)_k k!}
\frac{(m-k-y+2)_{k+y-1}(2x)_{k+y-1}}{(2x+m+1)_{k+y-1}(k+y-1)!}.
\tag\edh
\endalign
$$
Writing out the consecutive factors of the four Pochhammer symbols in the second fraction of the summand above, one obtains after some manipulation that it can be rewritten as
$$
(-1)^k
\frac{(m-y+2)_{y-1}(2x)_{y-1}}{(2x+m+1)_{y-1}(1)_{y-1}}
\frac{(-m+y-1)_k(2x+y-1)_k}{(2x+m+y)_k(y)_k}.
$$
Then, since $(y-k+1)_k=(-1)^k(-y)_k$,  (\edh) becomes
$$
\spreadlines{3\jot}
\align
&
\M(D_{x,y,y-1,m})=
\M(P_{y,y,x-1})\M(P_{m,m,x})
\frac{(m-y+2)_{y-1}(2x)_{y-1}}{(2x+m+1)_{y-1}(1)_{y-1}}
\\
&\ \ \ \ \ \ \ \ \ \ \ \ 
\times
\sum_{k=0}^{m-y+1}
\frac{(-y)_k(y+1)_k(-m+y-1)_k}{(y)_k(2x+m+y)_k k!}.
\tag\edi
\endalign
$$
Due to the presence of the third factor in the numerator of the summand in (\edi), the sum can be extended without change to infinity at the upper limit and can therefore be written as a hypergeometric function\footnote{ The hypergeometric function of parameters
$a_1,\dotsc,a_p$ and $b_1,\dotsc,b_q$ is defined by
$$
{}_p F_q\!\left[\matrix a_1,\dotsc,a_p\\ b_1,\dotsc,b_q\endmatrix;
z\right]=\sum _{k=0} ^{\infty}\frac {(a_1)_k\cdots(a_p)_k}
{k!\,(b_1)_k\cdots(b_q)_k} z^k\ .$$}
$$
\sum_{k=0}^{\infty}
\frac{(-y)_k(y+1)_k(-m+y-1)_k}{(y)_k(2x+m+y)_k k!}
=
{}_3F_2\!\left[\matrix{-y,\,y+1,\,-m+y-1}\\{y,\,2x+m+y}\endmatrix;1\right].
\tag\edj
$$
It turns out that the hypergeometric function in (\edj) can be evaluated in closed form using the classical transformation
$$
{}_3F_2\!\left[\matrix{-n,\,a,\,b}\\{c,\,d}\endmatrix;1\right]
=
\frac{(c+d-a-b)_n}{(c)_n}
{}_3F_2\!\left[\matrix{-n,\,d-a,\,d-b}\\{d,\,c+d-a-b}\endmatrix;1\right]
\tag\edk
$$
(see e.g. \cite{\PBM,Eq.\,83,\,p.\,539}).

Indeed, applying (\edk) with $n=y$, $a=y+1$, $b=-m+y-1$, $c=2x+m+y$ and $d=y$, we obtain
$$
{}_3F_2\!\left[\matrix{-y,\,y+1,\,-m+y-1}\\{y,\,2x+m+y}\endmatrix;1\right]
=
\frac{(2x+2m)_y}{(2x+m+y)_y}
{}_3F_2\!\left[\matrix{-y,\,-1,\,m+1}\\{y,\,2x+2m}\endmatrix;1\right].
\tag\edl
$$
But due to the numerator parameter $-1$ on the right hand side of (\edl), only the first two terms in the expansion of that hypergeometric function are non-zero, and we get
$$
\spreadlines{3\jot}
\align
{}_3F_2\!\left[\matrix{-y,\,-1,\,m+1}\\{y,\,2x+2m}\endmatrix;1\right]
=
\sum_{k=0}^{1}\frac{(-y)_k(-1)_k(m+1)_k}{(y)_k(2x+2m)_k}
&=
1+\frac{(-y)(-1)(m+1)}{y(2x+2m)}
\\
&=
\frac{2x+3m+1}{2x+2m}.
\tag\edm
\endalign
$$
By (\edi), (\edj), (\edl) and (\edm) we obtain
$$
\spreadlines{3\jot}
\align
&
\M(D_{x,y,y-1,m})=
\\
&\ \ \ \ \ \ 
\M(P_{y,y,x-1})\M(P_{m,m,x})
\frac{(m-y+2)_{y-1}(2x)_{y-1}}{(2x+m+1)_{y-1}(1)_{y-1}}
\frac{(2x+2m)_y}{(2x+m+y)_y}
\frac{2x+3m+1}{2x+2m}.
\tag\edn
\endalign
$$
Using the formula given in Corollary {\tdab}, it is routine to check that the expression on the right hand side of (\edn) agrees with the $z=y-1$ specialization of the product formula in Theorem~{\tbc}.
\epf

\topinsert
\centerline{\mypic{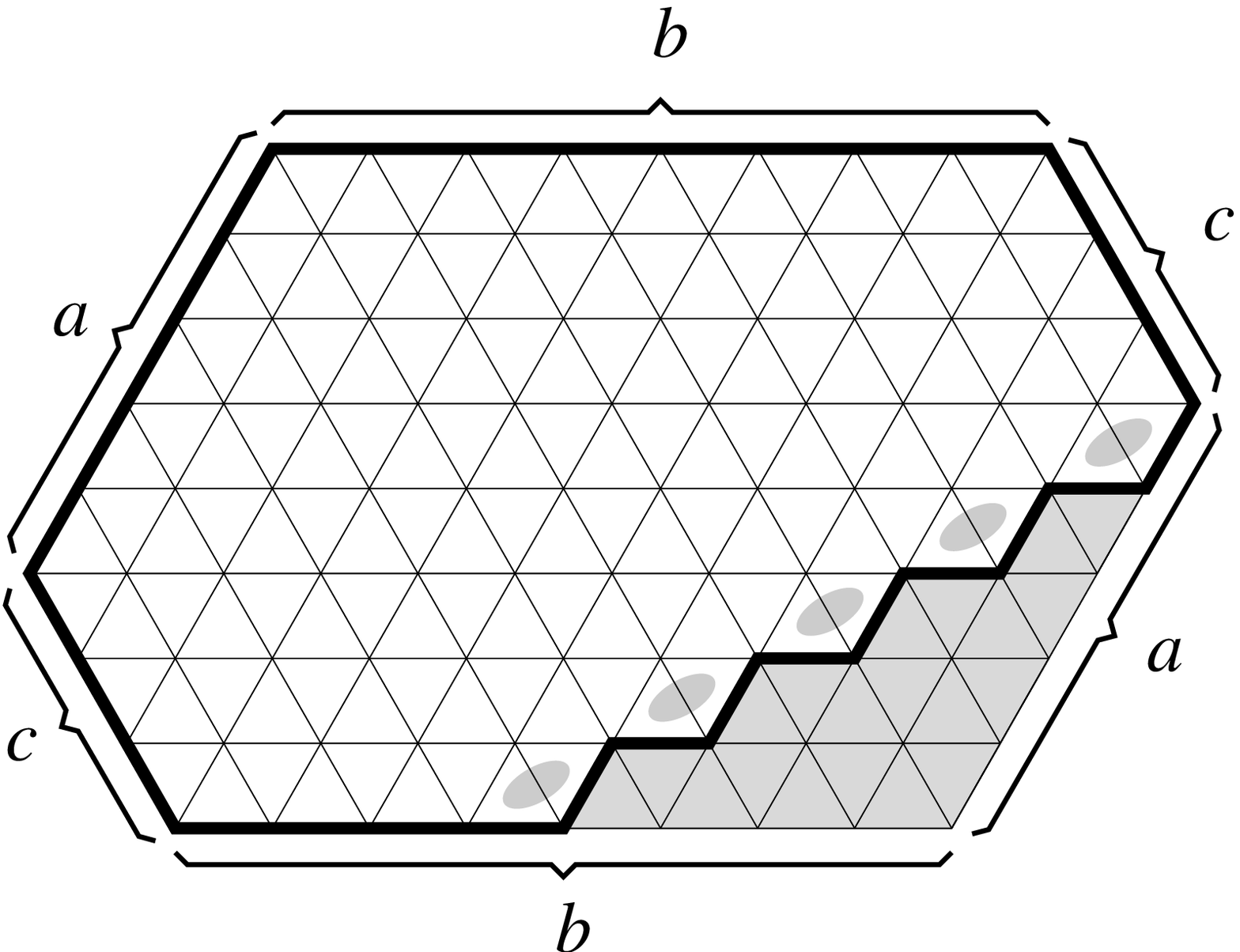}}
\medskip
\centerline{{\smc Figure~{\fde}. {\rm The region $P'_{a,b,c}$ for $a=5$, $b=8$, $c=3$.}}}
\endinsert

The proof of the special case $z=y-1$ of Theorem {\tbd} uses two variations of Corollary {\tdab} and Lemma {\tdb}, which we present below.



Let $P'(a,b,c)$ be the region obtained from $P(a,b,c)$ by weighting by $1/2$ the lozenge positions indicated in Figure {\fde} by shaded ellipses.

 \proclaim{Corollary {\tdab'}} For any non-negative integers $a$ and $c$ we have
$$
\M(P'_{a,a,c})=
\frac{(c+1/2)_a}{(2c)_a}
\prod_{1\leq i\leq j\leq a}
\frac{2c+i+j-2}{i+j-1}.
\tag\edo
$$

\endproclaim

\topinsert
\centerline{\mypic{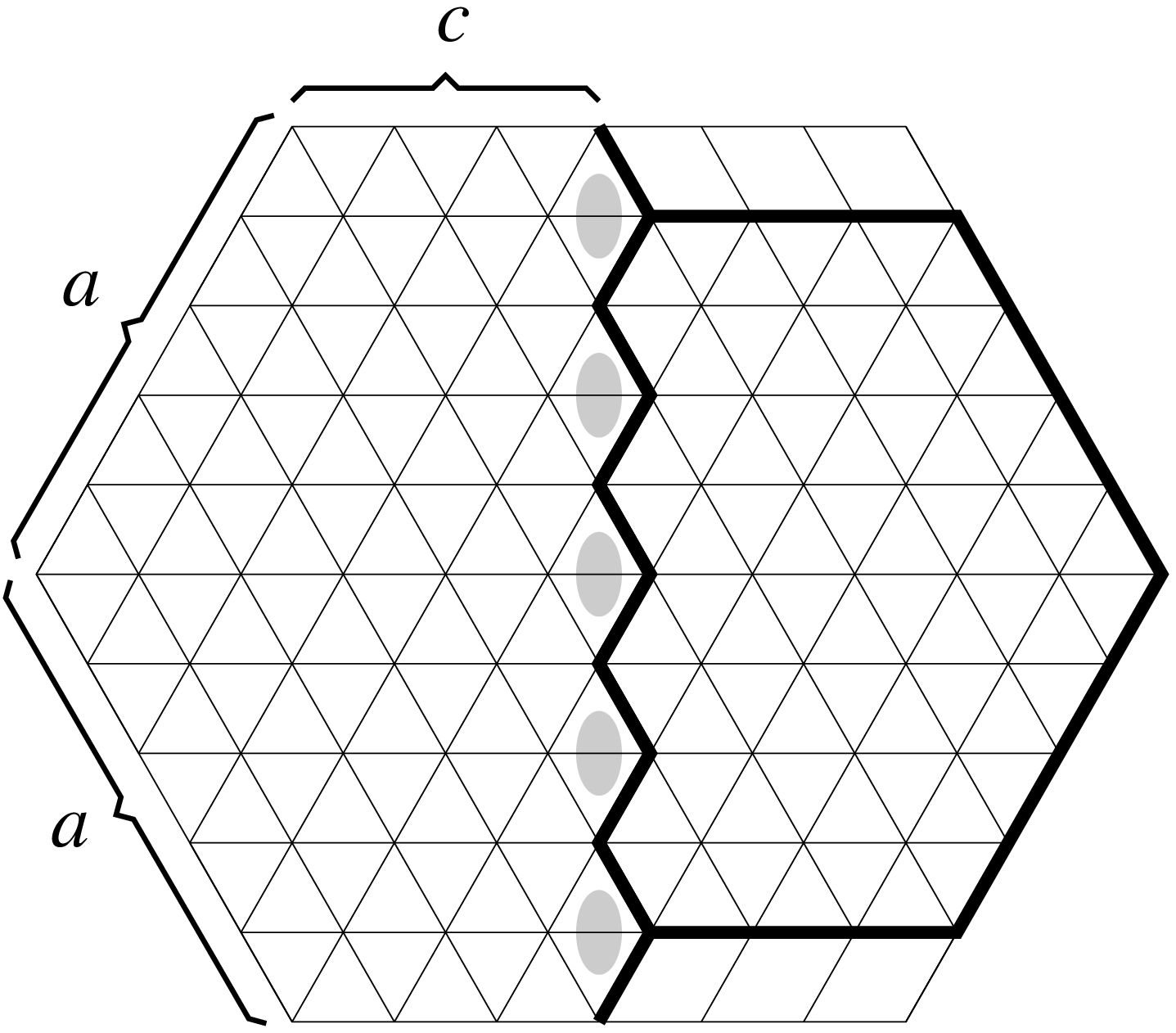}}
\medskip
\centerline{Figure~{\fdf}{. \rm $H_{a,a,c}$ and $P'_{a,a,c}$ for $a=5$, $c=3$.}}
\endinsert

\pf Apply the factorization theorem of \cite{\FT} to the dual graph of the hexagonal region $H_{a,2c,a}$ of side-lengths $a$, $2c$, $a$, $a$, $2c$, $a$ (clockwise from the northwestern side), with respect to its vertical symmetry axis (see Figure {\fdf}). After removing the forced lozenges from one of the resulting regions, we obtain
$$
\M(H_{a,2c,a})=2^a\M(P'_{a,a,c})\M(P_{a-1,a-1,c}).\tag\edp
$$
Using MacMahon's classical formula \cite{\MacM}
$$
\M(H_{a,b,c})=\prod_{i=1}^a\prod_{j=1}^b\prod_{k=1}^c\frac{i+j+k-1}{i+j+k-2}
$$
and Corollary {\tdab}, (\edp) yields a product expression for $\M(P'_{a,a,c})$, which is readily verified to agree with the right hand side of (\edo). \epf

\topinsert
\twoline{\mypic{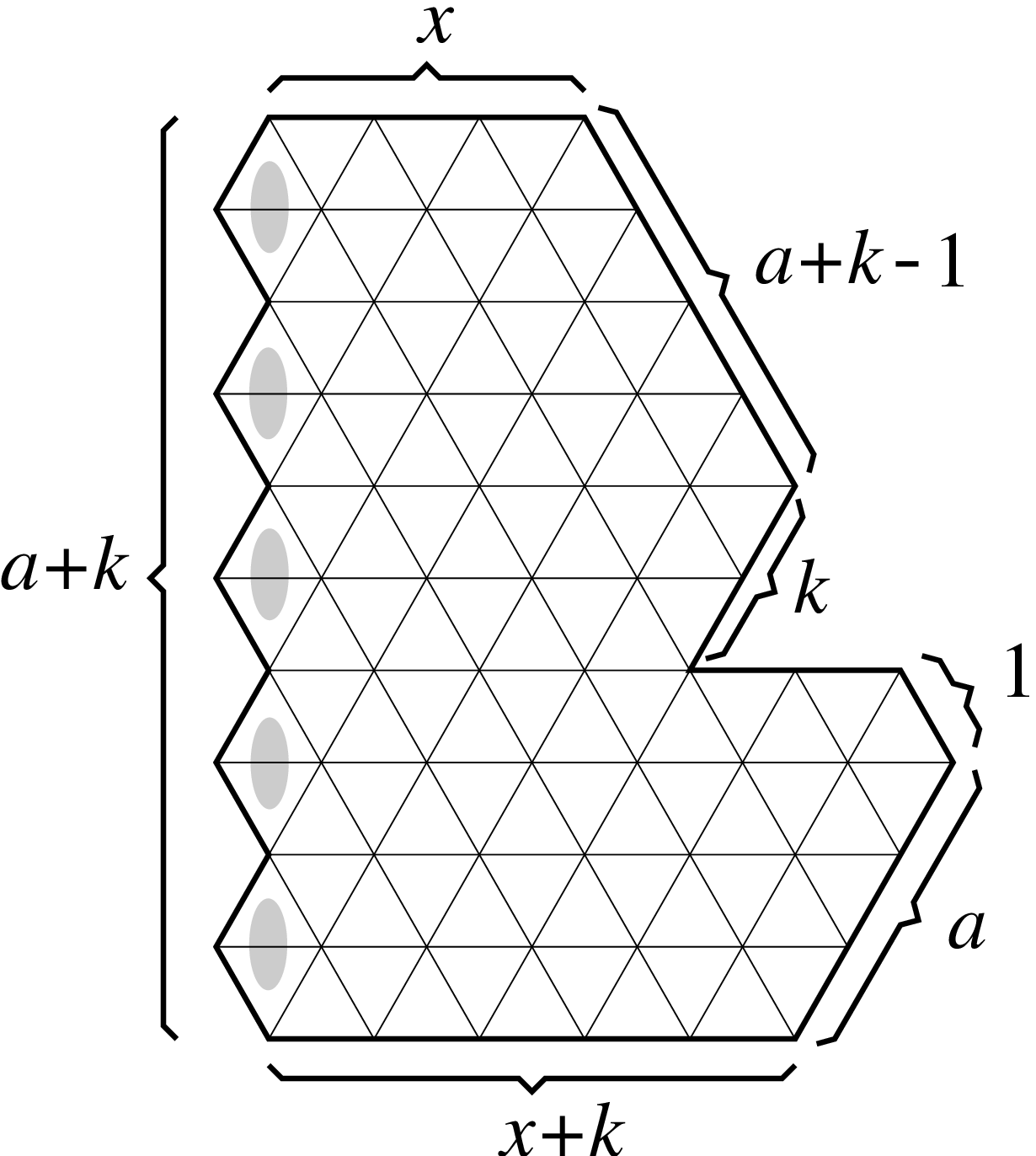}}{\mypic{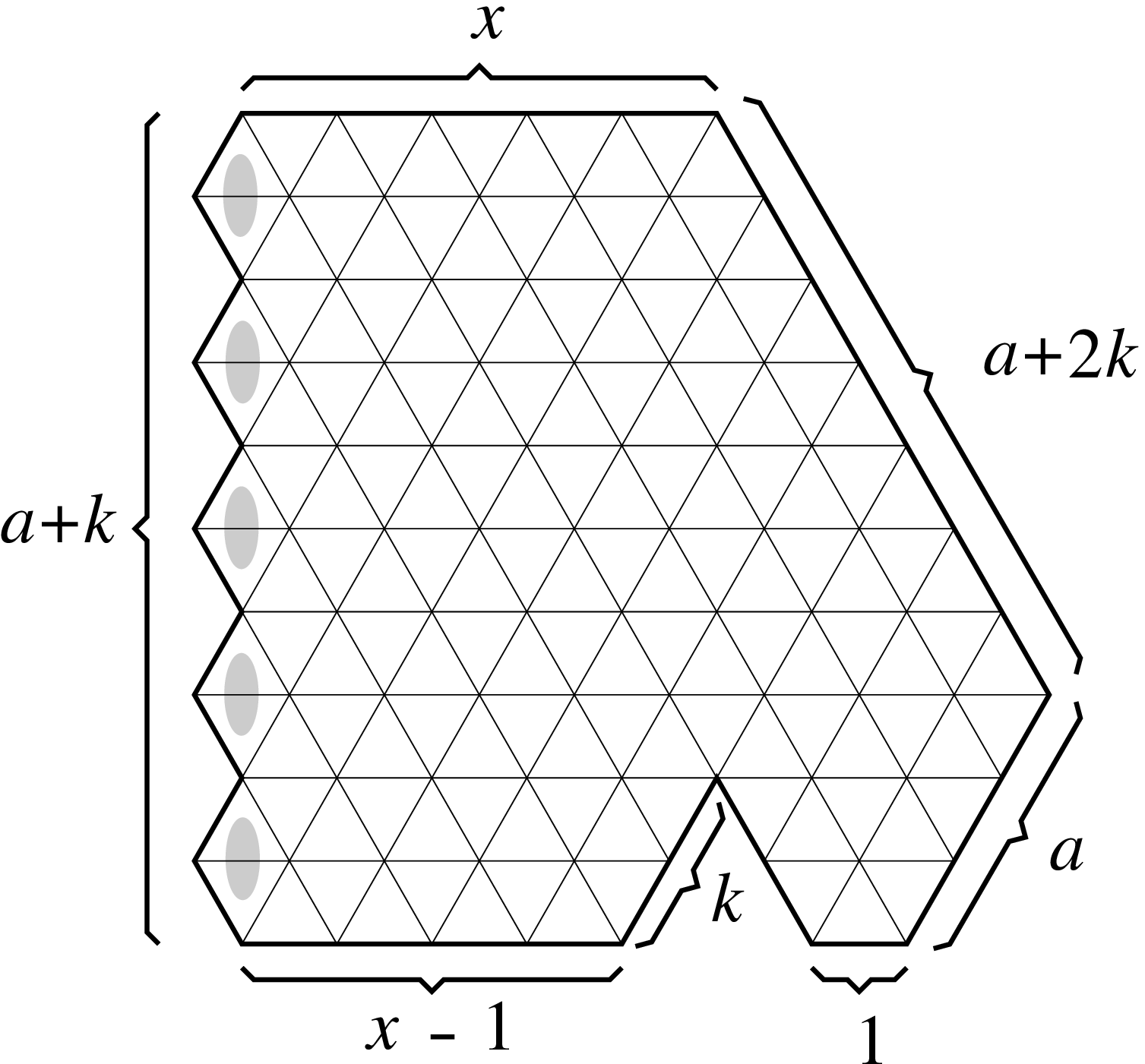}}
\medskip
\twoline{Figure~{\fdg}{. \rm $R'_{x,a,k}$ for $x=3$, $a=3$, $k=2$.}}
{Figure~{\fdh}{. $G'_{x,a,k}$.}}
\endinsert

In the same spirit, let $R'_{x,a,k}$ and $G'_{x,a,k}$ be the regions obtained from $R_{x,a,k}$ and $G_{x,a,k}$ by weighting by $1/2$ the tile positions indicated in Figures {\fdg} and {\fdh}. The counterpart of Lemma {\tdb} we need is the following.




\proclaim{Lemma {\tdb'}} For any non-negative integers $x$, $a$ and $k$ we have

$(${\rm a}$)$.
$$
\M(R'_{x,a,k})=
\frac{(a+1)_k (a+k+1)_k}{(2x+a+k)_k k!}\M(P'_{a+k,a+k,x}).
\tag\edq
$$

$(${\rm b}$)$.
$$
\M(G'_{x,a,k})=
\frac{(a+1)_k (2x-1)_k}{(2x+a+k)_k k!}\M(P'_{a+k,a+k,x}).
\tag\edr
$$

\endproclaim

\pf (a). The proof is perfectly analogous to the proof of Lemma {\tdb}(a). The base cases hold by Corollary {\tdab'}.

(b). This part can also be proved by the approach in the proof of Lemma {\tdb}(a). Alternatively, it follows directly from a classical result essentially due to Gelfand and Tsetlin (see e.g. \cite{\CLP, Proposition\,2.1}) and Lemma {\tdb}(b), using the factorization theorem of \cite{\FT}, as in the above proof of Corollary~{\tdab'}. \epf

{\it Proof of the $z=y-1$ case of Theorem {\tbd}.} The same arguments as in the proof of the special case $z=y-1$ of Theorem {\tbc} lead in this case to the following counterpart of (\edg):
$$
\M(D'_{x,y,y-1,m})=\sum_{k=0}^{m-y+1}\M(R'_{x-1,y-k,k})\M(G'_{x,m-(k+(y-1)),k+(y-1)}).
\tag\edes
$$
Since the formulas for $\M(R'_{x,a,k})$ and $\M(G'_{x,a,k})$ are obtained from the formulas for $\M(R_{x,a,k})$ and $\M(G_{x,a,k})$ by simply replacing $x$ by $x-1/2$ (compare the formulas in Lemmas {\tdb} and {\tdb'}), and also the formula in Corollary {\tdab'} is obtained from the formula in Corollary {\tdab} by replacing $c$ by $c-\frac12$, it follows by (\edes) and (\edg) that $\M(D'_{x,y,y-1,m})$ is equal to the expression obtained from the $z=y-1$ specialization of the right hand side of (\ebc). Since replacing $x$ by $x-1/2$ in the expression on the right hand side of (\ebc) transforms it into the expression on the right hand side of (\ebd) (in fact, this is true in general, and not just in the $z=y-1$ case), the proof of the special case $z=y-1$ of Theorem {\tbd} is completed. \epf

\mysec{5. The proofs of Theorems {\tbc} and {\tbd}}

Our proofs are inductive, and rely on the recurrences stated in Lemmas {\tca} and {\tcb}.

{\it Proof of Theorem {\tbc}.} We prove (\ebc) by induction on $x$. By (\ecg), for any positive integers $x,y,z,m$ with $z<y$, $z<m$ and $y-z\leq x$, we have 
$$
\spreadlines{3\jot}
\align
\M(D_{x+1,y,z-1,m-1})=
\frac{\M(D_{x,y,z,m})\M(D_{x,y-1,z-1,m-1})-\M(D_{x,y-1,z-1,m})\M(D_{x,y,z,m-1})}
{\M(D_{x-1,y-1,z,m})}.
\tag\eea
\endalign
$$

Equivalently, for any non-negative integers $x,y,z,m$ with $x\geq2$, $y\geq1$, $z<y-1$, $z<m$ and $y-z\leq x$, we have 
$$
\spreadlines{3\jot}
\align
\M(D_{x,y,z,m})
\!=\!
\frac{\M(D_{x-1,y,z+1,m+1})\M(D_{x-1,y-1,z,m})
\!-\!
\M(D_{x-1,y-1,z,m+1})\M(D_{x-1,y,z+1,m})}
{\M(D_{x-2,y-1,z+1,m+1})}.
\tag\eeb
\endalign
$$
Thus the cases when (\eeb) cannot be applied are 
\medskip
\ \ $(i)$ $x=0$

\ $(ii)$ $x=1$

$(iii)$ $y=0$

\ $(iv)$ $z=y-1$

\ \ $(v)$ $z=y$

\ $(vi)$ $z=m$
\medskip
\flushpar
(Note that the quantity $x-(y-z)$, which recall is the depth of the notch at the top of our region, is the same across all six terms in (\eeb), so provided it is non-negative for the region on the left, it is non-negative also for all five regions on the right.)

Cases $(i)$ and $(ii)$ are the base cases of our induction. We verify cases $(iii)-(vi)$ by separate arguments. 

$(i)$. If $x=0$, then $0\leq x-(y-z)=z-y$, and since we always have $z\leq y$, it follows that $z=y$. So this case follows from case $(v)$ below.

$(ii)$. $x=1$ implies $1=x\geq x-(y-z)\geq 0$, so the depth $x-y+z$ of the notch at the top of our region is either 0 or 1. If it is 0, then we have $y-z=1$, so $z=y-1$, and the statement follows by Section 4. On the other hand, if $x-y+z=1$, then it follows that $z=y$, which reduces this case to case $(v)$.

\topinsert
\centerline{\mypic{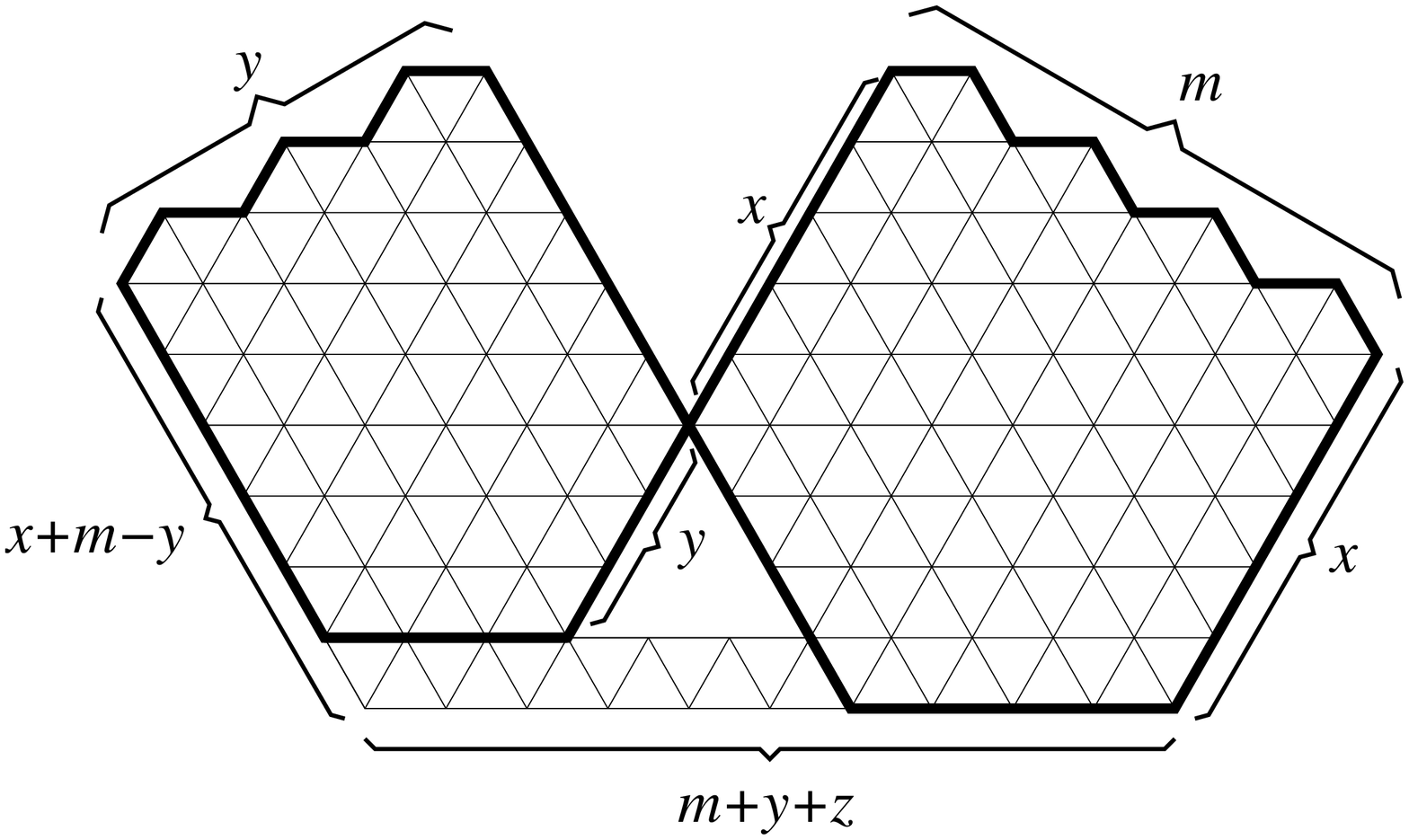}}
\medskip
\centerline{Figure~{\fea}{. \rm $D_{5,3,3,4}$.}}
\endinsert

$(iii)$. If $y=0$, then necessarily $z=0$ as well, so this follows from case $(v)$.

$(iv)$. This case was proved in Section 4.

$(v)$. If $y=z$, then $D_{x,y,z,m}$ is of the type illustrated in Figure {\fea}. Consider an arbitrary tiling of $D_{x,y,y,m}$, and consider its encoding as a family of $x+y$ non-intersecting paths of lozenges, $x$ of which start on the southeastern side, and the remaining $y$ starting on the northwestern side of the triangular gap. The latter must end at the $y$ southeast facing unit segments along the top left zig-zag boundary of $D_{x,y,y,m}$, and therefore the former must end at the $x$ unit segments along the southeastern side of the notch at the top of $D_{x,y,y,m}$. It follows that the two regions inside the thick contours in Figure {\fea} are tiled internally (i.e., no lozenge can cut across their boundaries). Since the remaining portion of $D_{x,y,y,m}$ is a rhombus, and therefore uniquely tileable, it follows then $\M(D_{x,y,y,m})$ is equal to the product of the number of lozenge tilings of those two regions. However, both these regions are of the type covered by Corollary {\tdab}. Their precise parameters can be read from Figure {\fea}, and we obtain
$$
\M(D_{x,y,y,m})=\M(P_{y,y,x})\M(P_{m,m,x}).
\tag\eec
$$
Expressing the factors on the right hand side above by Corollary {\tdab}, (\eec) gives a product formula for $\M(D_{x,y,y,m})$. It is straightforward to check that this agrees with the $z=y$ specialization of formula (\ebc).

$(vi)$. This follows by interchanging the roles of $y$ and $m$ in case $(v)$.

\medskip

For the induction step, suppose (\ebc) holds for all the $D$-regions with $x$-parameter strictly less than $x$, and consider the region $D(x,y,z,m)$. If any of the equations $(i)-(vi)$ above hold, (\ebc) follows by the above arguments. Otherwise we have $x\geq2$, $y\geq1$, $z<y-1$, $z<m$ and $y-z\leq x$, and thus (\eeb) holds. Using the induction hypothesis, all five $D$-regions on the right hand side of (\eeb) have the number of their lozenge tilings given by formula (\ebc). To complete the induction step it suffices to show that the expression on the right hand side of (\ebc) also satisfies recurrence (\eeb). However, this follows by Lemma {\tfa}. This completes the induction step, and therefore the proof of Theorem {\tbc}. \epf

{\it Proof of Theorem {\tbd}.} By Lemma {\tcb}, the numbers $\M(D'_{x,y,z,m})$ also satisfy the recurrence (\eeb). Just as in the proof of Theorem {\tbc}, this recurrence can be applied unless one of the equalities $(i)-(vi)$ above holds. The verification of these special cases works exactly as in the proof of Theorem {\tbc}, with the one change that instead of Corollary {\tdab} we use Corollary {\tdab'} (note that case $(iv)$ still follows by Section 4, as that section contains the proof of the special case $z=y-1$ of both Theorems {\tbc} and {\tbd}).

The induction step follows precisely as in the proof of Theorem {\tbc}, using the fact that the expression on the right hand side of (\ebd) satisfies recurrence (\eeb) by Lemma {\tfa}. \epf

\mysec{6. Verification that expressions (\ebc) and (\ebd) satisfy recurrence (\ecg)}

Let us denote by $f(x,y,z,m)$ the expression on the right hand side of (\ebc), and by $g(x,y,z,m)$ the expression on the right hand side of (\ebd), with the convention (as in the statements of Theorems {\tbc} and {\tbd}) that if the upper index in a product is larger than the lower index, then that product is taken to be 1.

The object of this section is to prove the following result.

\proclaim{Lemma \tfa} Both $f(x,y,z,m)$ and $g(x,y,z,m)$ satisfy recurrence $(\ecg)$, i.e., for any non-negative integers $x,y,z$ and $m$ with $z < y$ and $z < m$ we have
$$
\spreadlines{3\jot}
\align
f(x,y,z,m)\, f(x,y-1,z-1,m-1) 
=
&f(x,y-1,z-1,m)\, f(x,y,z,m-1)
\\
 + 
&
f(x+1,y,z-1,m-1)\, f(x-1,y-1,z,m)
\tag\efa
\endalign
$$
and
$$
\spreadlines{3\jot}
\align
g(x,y,z,m)\, g(x,y-1,z-1,m-1) 
=
&\,
g(x,y-1,z-1,m)\, g(x,y,z,m-1)
\\ 
+ 
&g(x+1,y,z-1,m-1)\, g(x-1,y-1,z,m).
\tag\efb
\endalign
$$

\endproclaim

\pf Since, as noted in Section 1 (see Remark 1), $g(x,y,z,m)=f(x-\frac12,y,z,m)$, it suffices to verify (\efa).

By dividing both sides by its left hand side, (\efa) can be rewritten equivalently as
$$
\frac
{\dfrac{f(x,y-1,z-1,m)}{f(x,y-1,z-1,m-1)}}
{\dfrac{f(x,y,z,m)}{f(x,y,z,m-1)}}
+
\frac
{\dfrac{f(x+1,y,z-1,m-1)}{f(x,y-1,z-1,m-1)}}
{\dfrac{f(x,y,z,m)}{f(x-1,y-1,z,m)}}
=1.
\tag\efc
$$
Set
$$
r_1(x,y,z,m):=\frac{f(x,y,z,m)}{f(x,y,z,m-1)}
\tag\efd
$$
and
$$
r_2(x,y,z,m):=\frac{f(x,y,z,m)}{f(x-1,y-1,z,m)}.
\tag\efe
$$
Then (\efc) becomes
$$
\frac{r_1(x,y-1,z-1,m)}{r_1(x,y,z,m)}
+
\frac{r_2(x+1,y,z-1,m-1)}{r_2(x,y,z,m)}
=1.
\tag\eff
$$
We divide our analysis in three cases, according to the residues of $m-z$ modulo 3. If $m-z$ is a multiple of 3, we obtain directly from (\efd) and (\ebc) that
$$
\spreadlines{3\jot}
\align
&
r_1(x,y,z,m)
=
\frac{m!(2z-2y+2z+2m)}{(2m)!}
\frac{(2x+m+3z+2(m-z-1)-y+3)_{y-m}}{ (y+(m-z-1)-z+1)_{z}}
\\
&\ \ \ \ \ \
\times
\frac{1}{(2x-2y+5z+6((m-z)/3-1)+6)}
\frac{(2x-2y+3z+(m-z-1)+2)_{y-z}}{((m-z-1)+1)_z}
\\
&\ \ \ \ \ \ 
\times
\prod_{i=0}^{m-2}
\frac{(2x-2y+2z+2i+2i)_{m-i}}{(2x-2y+2z+2i+2)_{m-1-i}}
\prod_{i=0}^{m-z-2}
\frac{(2x+m+3z+2i-y+3)_{y-m}}{(2x+(m-1)+3z+2i-y+3)_{y-(m-1)}}
\\
&\ \ \ \ \ \
\times
\prod_{i=0}^{(m-z)/3-1}
\frac{(2x+y+2z+3i+3)_{3m-3z-2-9i}}{(2x+y+2z+3i+3)_{3(m-1)-3z-2-9i}}
\\
&\ \ \ \ \ \
\times
\prod_{i=0}^{(m-z)/3-1}
\frac{(2x-2y+5z+6i+5)_{3m-3z-5-9i}}{(2x-2y+5z+6i+5)_{3(m-1)-3z-5-9i}}
\\
&\ \ \ \ \ \
\times
\prod_{i=0}^{(m-z)/3-1}
\frac{(2x+y-z+3(m-1)-6i-1)}{(2x+y-z+3(m)-6i-1)}.
\tag\efg
\endalign
$$
%
%
When replacing $y$ by $y-1$ and $z$ by $z-1$ in (\efg) and dividing by the expression (\efg), most factors simplify out. All that remains is
$$
\frac{r_1(x,y-1,z-1,m)}{r_1(x,y,z,m)}
=
\frac
{(m-z)(2 x + y + z + m)(2 x - 2 y + 3 z + 2 m)}
{(y + m - 2 z) (2x+z+2m)(2x-y+3z+m)}.
\tag\efh
$$
%

Remarkably, despite the minor variations (\efg) takes on when $m-z$ is congruent to 1 or 2 modulo 3, equation (\efh) turns out to hold without change in both these remaining cases. Thus (\efh) holds in general. 

A similar analysis leads to the fact that
$$
\frac{r_2(x+1,y,z-1,m-1)}{r_2(x,y,z,m)}
=
\frac
{(y - z)(2 x + 3 z)(2x-y+z+3m)}
{(y + m - 2 z)(2x+z+2m)(2x-y+3z+m)}
\tag\efi
$$
for all residues of $m-z$ modulo 3.
%

Since the sum of the two fractions on the right hand sides of (\efh) and (\efi) is readily seen to equal 1, (\eff) follows by (\efh) and (\efi), thus proving (\efa). \epf

\mysec{7. Concluding remarks}

We have seen in this paper another instance of the power of Kuo's graphical condensation method, two other related applications of which were presented in \cite{\anglepap} and \cite{\ffp}. Our approach also illustrates the well-known principle that sometimes it is easier to prove a generalization of a conjecture than the conjecture itself --- indeed, we do not have a proof of the original Conjectures {\tba} and {\tbb} that proves them directly, without first generalizing them. In particular, the graphical condensation method only works for the more general families of regions of Theorems {\tbc} and {\tbd}.

We end this paper by emphasizing the curious fact that $\M(D'_{x,y,z,m})$ is equal to what results from the product expression for $\M(D_{x,y,z,m})$ when replacing $x$ by $x-\frac12$ (see Remark~1 at the end of Section 2). It would be interesting to have a direct justification of this.

\mysec{References}
{\openup 1\jot \frenchspacing\raggedbottom
\roster

\myref{\And}
  G. E. Andrews, Plane partitions (III): The weak Macdonald
conjecture, {\it Invent. Math.} {\bf 53} (1979), 193--225.

\myref{\FT}
  M. Ciucu, Enumeration of perfect matchings in graphs with reflective symmetry,
{\it J. Combin. Theory Ser. A} {\bf 77} (1997), 67--97.

\myref{\cutcor}
  M. Ciucu and C. Krattenthaler, Enumeration of lozenge tilings of hexagons with cut off corners, {\it J. Combin. Theory Ser. A} {\bf 100} (2002), 201-231 .

\myref{\anglepap}
  M. Ciucu and I. Fischer, A triangular gap of size two in a sea of dimers on a
$60^\circ$ angle, {\it J. Phys. A: Math. Theor.} {\bf 45}  (2012), 494011.

\myref{\ffp}
  M. Ciucu and C. Krattenthaler, A dual of MacMahon's theorem on plane partitions, {\it Proc. Natl. Acad. Sci. USA}  {\bf 110} (2013), 4518--4523.

\myref{\CLP}
  H. Cohn, M. Larsen, and J. Propp, The shape of a typical boxed plane partition,{\it New York J. of Math.} {\bf 4} (1998), 137--165.

\myref{\Kuo}
Kuo, Eric H. Applications of graphical condensation for enumerating matchings and tilings. Theoret. Comput. Sci. 319 (2004), no. 1-3, 29–57.

\myref{\Kup}
  G. Kuperberg, Symmetries of plane partitions and the permanent-de\-ter\-mi\-nant meth\-od, {\it J. Combin. Theory Ser. A} {\bf 68} (1994), 115--151.




\myref{\Tri} 
 T. Lai, Quartered Aztec rectangles with holes and quartered hexagons, preprint, 2013.

\myref{\MacM}
  P. A. MacMahon, Memoir on the theory of the partition of numbers---Part V. Partitions in two-dimensional space, {\it Phil. Trans. R. S.}, 1911, A.

\myref{\Proctor}
  R. A. Proctor, Odd symplectic groups, {\it Invent. Math.} {\bf 92} (1988),
307--332.

\myref{\PBM}
  A.P. Prudnikov, Yu.A. Brychkov and O.I. Marichev, ``Integrals and series,'' vol. 3,
Gordon and Breach Science Publishers, New York, 1986.

\myref{\Sta}
  R. P. Stanley, Ordered structures and partitions, {\it Memoirs of the Amer. Math. Soc.,} 
no. 119 (1972).

\myref{\Ste}
  J. R. Stembridge, Nonintersecting paths, Pfaffians and plane partitions,
{\it Adv. in Math.} {\bf 83} (1990), 96--131.

\endroster\par}

\enddocument